\documentclass{lqr-preprint}
\newcommand{\preprintbib}{references}
\jmlryear{2026}
\ifdefined\pdfminorversion\pdfminorversion=4\fi
\usepackage{times}
\usepackage{mathtools}
\usepackage{enumitem}
\usepackage{microtype}
\usepackage{booktabs}

\hypersetup{hidelinks}
\newtheorem{assumption}{Assumption}
\let\lqroriginalproof\proof
\let\endlqroriginalproof\endproof
\renewenvironment{proof}[1][Proof]{%
  \renewcommand{\proofname}{#1}\lqroriginalproof
}{\endlqroriginalproof}

\usepackage{wrapfig}
\makeatletter
\newcommand{\narrowcaptionstyle}{%
  \long\def\@makecaption##1##2{\vskip 1pt{\small\parskip 0pt ##1: ##2\par}}}
\makeatother

\newcommand{\dt}{{\Delta t}}
\newcommand{\E}{\mathbb{E}}

\title[Continuous-Time LQR from a Single Off-Policy Batch]{Continuous-Time LQR from a Single Discrete-Time Off-Policy Batch}
\author[Yao and Zhu]{%
\Name{Chaorui Yao\nametag{\textsuperscript{1}}} \Email{chaorui@ucla.edu}\\
\Name{Yuhua Zhu\nametag{\textsuperscript{2}\thanks{Corresponding author.}}} \Email{yuhuazhu@ucla.edu}\\
\addr \textsuperscript{1}Department of Electrical and Computer Engineering, University of California, Los Angeles\\
\addr \textsuperscript{2}Department of Statistics and Data Science, University of California, Los Angeles}

\begin{document}
\maketitle
\begin{abstract}
We study discounted continuous-time linear--quadratic
regulation (LQR) with unknown drift parameters and additive noise, using
a single batch of discrete-time off-policy transitions.
Discrete-time methods target a policy tied to the
training interval, while shortening this interval can make learning from
finite data more sensitive to noise. Learning an accurate continuous-time
controller therefore requires balancing discretization and finite-sample
errors while preserving stability throughout policy improvement.
We develop
a drift-only policy iteration algorithm, termed PhiBE-PI, using one
quadratic state--action critic for policy evaluation and explicit greedy
improvement. Reusing
the batch makes later policies depend on the data. Exploiting the
quadratic structure of LQR, we show that the data enter every policy
evaluation only through one policy-independent empirical matrix, so a
single high-probability bound controls all iterates.

For sufficiently large discount~$\beta$, small observation
intervals, and a sufficiently large batch, we prove that
all iterates are stabilizing for the discounted problem and converge
quadratically to a data-dependent fixed point.
The gain bound separates iteration, finite-sample, and discretization errors.
In this regime, the finite-sample and discretization
terms scale as $1/(\beta\sqrt{n\Delta t})$ and $\Delta t/\beta$,
respectively. The bounds
thus reveal a tradeoff: at a fixed batch size, finer observations reduce
discretization error but increase finite-sample error.
Because PhiBE-PI approximates the continuous-time optimal
gain rather than a discrete-time one, its policy can remain near-optimal
when deployed at faster control rates. Experiments confirm this advantage
over discrete-time Bellman-equation baselines when the observation
interval is small relative to the system's time scale.
\end{abstract}

\section{Introduction}
\label{sec:intro}
Continuous-time linear quadratic regulation (LQR) approximates
nonlinear optimal control problems near equilibria through linearized dynamics
and quadratic costs. This motivates its use in stabilization and tracking tasks,
including robotic balancing and quadrotor control \citep{TedrakeUnderactuated}.
When the dynamics are unknown, discrete observation creates
two challenges. First, the controller must be learned from a finite batch of
discrete observations while targeting the continuous-time optimum. Second, its
feedback law must be implemented through discrete updates: frequent updates
enable prompt responses to state changes, but their frequency is limited by
the observation rate.

\begin{wrapfigure}{r}{0.47\textwidth}
\vspace{-\intextsep}
\narrowcaptionstyle
\centering
\includegraphics[width=\linewidth]{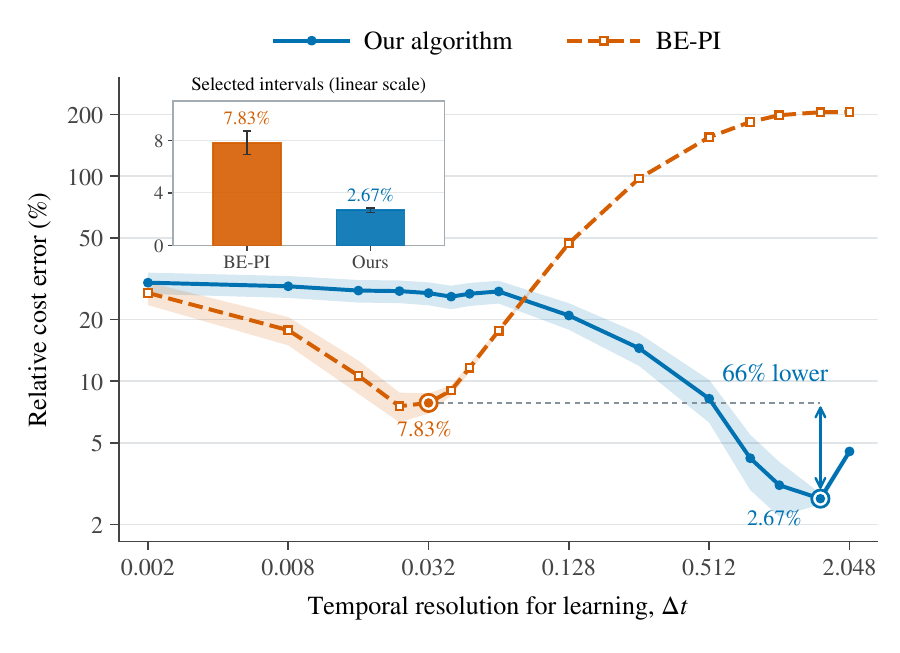}\par
\caption{Relative cost error versus the training
interval $\Delta t$ for a scalar LQR system. Both methods use the same data
and deploy at the observation interval $\delta=0.002$. Curves and bands show
means and 95\% confidence intervals over 1000 paired trials. Rings and the
inset compare intervals selected before validation. Failed fits retain the
common initial controller; all trials are included. BE-PI is discrete-time
Bellman-equation policy iteration. Details are given in
Appendix~\ref{app:motivation-setup}.}
\label{fig:learning-temporal-resolution}
\end{wrapfigure}

Even when observations are collected at the highest available
rate, learning at this frequency can be challenging for standard reinforcement learning
methods \citep{TallecBlierOllivier2019,ParkKimKim2021SafeActionRepetition}.
Provided the control input is held constant between
them, observations farther apart can be paired for training while the
controller is deployed at the original observation rate.
The temporal resolution used for learning
is therefore a design choice. As illustrated in
Figure~\ref{fig:learning-temporal-resolution}, at separately selected
training resolutions, our algorithm achieves a relative cost error of
$2.67\%$, compared with $7.83\%$ for BE-PI, a reduction of approximately
$66\%$. This motivates studying how temporal resolution
and finite data jointly affect the accuracy of learned continuous-time LQR
controllers.

Two routes are available. The discrete-time route treats each
training transition as one step of a Markov decision process (MDP) and applies
established methods such as least-squares policy iteration
\citep{KrauthTuRecht2019}. The optimization problem this defines is tied
to the training interval: one step of the MDP is one interval with the action
held fixed, so its solution is the best controller that updates at exactly
that rate. This solution generally changes with the interval and need not
remain optimal when deployed at a different frequency. Shrinking the interval
does not remove this dependence in practice. The problem itself is well
behaved in this limit: under standard assumptions, its optimal gain approaches
the continuous-time counterpart as the interval tends to zero
\citep{BourdinTrelat2021UnifiedRiccati}. But this gain also becomes
harder to learn for methods that select actions by comparing action values,
such as Q-learning and least-squares policy iteration. An action held for a
single interval changes its action value only by an amount proportional to
the interval. A fixed level of estimation error in the action values
therefore has a growing effect on the selected action
\citep{Baird1994AdvantageUpdating,TallecBlierOllivier2019,JiaZhou2023}.

The continuous-time route instead targets the original continuous-time
objective. With known dynamics, Kleinman's policy iteration
preserves stability and converges to the optimal gain from a stabilizing
initialization, with a locally quadratic rate under standard assumptions
\citep{Kleinman1968}.
With unknown dynamics, model-based methods estimate system
matrices before synthesizing a controller. From discrete observations,
\citet{BaseiGuoHuZhang2022} incur discretization error in episodic
finite-horizon LQR, whereas \citet{ZhouLiZhang2025FiniteSample} bypass it in
online LQR at a fixed sampling interval chosen using prior knowledge of the
dynamics. In contrast, the model-free PhiBE
\citep{Zhu2024PhiBE,ZhuZhangZhang2025OptimalPhiBE}
approximation underlying our method
achieves zero discretization error in the optimal gain under certain conditions,
without identifying the dynamics.
Among model-free methods using fixed trajectory data,
\citet{LopezMuller2023CTLQR} prove stability and quadratic convergence, while
\citet{PangJiang2023Stationary} establish near-optimality for stochastic
systems from sufficiently long trajectories. Both analyses rely on
continuous-trajectory integrals without quantifying the error introduced by
discrete observations, so their guarantees do not establish lower
discretization error than the discrete-time MDP approach.

We study discounted continuous-time LQR with unknown drift parameters and
additive noise. The learner receives a single batch of one-step off-policy
transitions, each generated by holding an action constant over a fixed
training interval. Carrying out policy iteration in this
setting poses three technical challenges. First, policy iteration must evaluate
and improve state-dependent feedback policies using only endpoint observations
collected under constant actions. Second, finite-sample
errors can make policy evaluation or greedy
improvement ill defined, and reusing the batch makes subsequent policies
dependent on the same estimation errors. For discrete-time
LQR, \citet{KrauthTuRecht2019} use fresh data for each policy evaluation,
whereas \citet{AntosSzepesvariMunos2008} allow data
reuse but assume bounded rewards and finitely many actions, excluding LQR.
Finally, convergence alone does not ensure accuracy relative to
the continuous-time optimum: at a fixed batch size, shorter training intervals
reduce discretization bias but can increase finite-sample error.
Prior work quantifies discretization and finite-sample
errors for fixed-policy evaluation
\citep{ZhangEtAl2023TemporalResolution,Mou2025ContinuousTimePolicyEvaluation},
but does not cover repeated policy improvement. This raises the
following question: \textbf{Can we learn a continuous-time LQR
controller from a single batch of discrete observations while preserving
discounted stability and quadratic convergence, with explicit accuracy bounds
in terms of temporal resolution and sample size?}

We develop a model-free policy iteration algorithm for discounted continuous-time LQR.
We approximate the unknown drift using conditional means of normalized
state increments, following the moment-based construction of PhiBE
\citep{Zhu2024PhiBE}. Our main contributions are as follows.

\begin{itemize}[leftmargin=*,itemsep=3pt,parsep=0pt,topsep=4pt,partopsep=0pt]
\item \textbf{Quadratic state--action critic.} Building on
\citet{companion_ct_q}, we formulate a quadratic state--action critic for
continuous-time LQR. It supports off-policy learning from discrete
observations while preserving the action-dependent signal as the training
interval tends to zero.

\item \textbf{Policy iteration from a single batch.} We reuse
a single batch of noisy discrete off-policy transitions, evaluating each policy
through a linear system and improving it through an explicit greedy gain update.

\item \textbf{Stability and quadratic convergence.} For sufficiently large
discount, sufficiently small training intervals, and a sufficiently large
batch, we prove that all policy evaluations and improvements are well defined,
all iterates remain stabilizing for the discounted problem, and the gains
converge quadratically to a limiting gain determined by the observed batch.
All evaluations share one policy-independent empirical matrix;
controlling its estimation error and using a quadratic error identity yield
these guarantees on one high-probability event for all data-dependent iterates.

\item \textbf{End-to-end error bounds.} We bound the learned gain's distance
from the continuous-time optimum by separating iteration, finite-sample, and
discretization errors. In the stated regime, the leading finite-sample and
discretization terms scale as $1/(\beta\sqrt{n\Delta t})$ and
$\Delta t/\beta$,
respectively, where $n$ is the batch size, $\Delta t$ the training interval,
and $\beta$ the discount rate.

\item \textbf{Discretization--finite-sample tradeoff.}
Our error bounds characterize the tradeoff between discretization and
finite-sample errors. Numerical experiments show that choosing an
appropriate temporal resolution for learning allows our algorithm to
learn a more accurate controller from the same batch of high-frequency
observations, while keeping deployment at the original observation rate
(Figure~\ref{fig:learning-temporal-resolution}).
\end{itemize}

\section{Problem Formulation and Algorithm}
\label{sec:lqr-case-study}

\subsection{Continuous-time control from sampled data}
\label{subsec:lqr-1d-setup}

Let $X_t\in\mathbb R^d$ denote the state and $u_t\in\mathbb R^m$ the
control input. We consider the linear stochastic system and quadratic reward
\begin{equation}
\label{eq:md-lqr-dynamics}
\begin{gathered}
dX_t=(AX_t+Bu_t)\,dt+\Gamma\,dW_t,\\
\mathcal R(x,u)=-(x^\top Qx+u^\top Ru),\quad Q\succeq0,\ R\succ0.
\end{gathered}
\end{equation}
Here $A\in\mathbb R^{d\times d}$ and $B\in\mathbb R^{d\times m}$ are
the state and input matrices, $W_t$ is a standard $d_w$-dimensional
Brownian motion, and $\Gamma\in\mathbb R^{d\times d_w}$ is the diffusion
matrix. The cost weights $Q,R$ are symmetric.
For discount rate $\beta>0$, a feedback policy $u_t=-KX_t$ with
$K\in\mathbb R^{m\times d}$ has value
\begin{equation}
\label{eq:lqr-ct-objective}
V_K^{\mathrm{ct}}(x)=\E_x\!\left[\int_0^\infty e^{-\beta t}
\mathcal R(X_t,-KX_t)\,dt\right].
\end{equation}
The expectation $\E_x$ is conditional on $X_0=x$.
We call a feedback gain $K$ discounted-stabilizing if
$A-BK-\beta I/2$ is Hurwitz. This condition ensures that the discounted
value in \eqref{eq:lqr-ct-objective} is finite.

\paragraph{Data model.}
The drift coefficients $A,B$ are unknown. Our goal is to learn
the optimal continuous-time feedback gain $K^\star_{\mathrm{ct}}$.
We observe one off-policy batch $\mathcal D=\{(X_i,u_i,R_i,X_i')\}_{i=1}^n$,
where $(X_i,u_i)\overset{\mathrm{i.i.d.}}{\sim}\mu$ and
$R_i=\mathcal R(X_i,u_i)$.
Each $X_i'$ is generated independently by
system~\eqref{eq:md-lqr-dynamics}, starting from $X_i$ and holding $u_i$
constant over an interval of length $\Delta t$. We assume a discounted-stabilizing initial
gain $K_0$, independent of the batch.

\subsection{Drift-only Q-function and its approximation}
\label{subsec:lqr-q-approximation}

We formulate drift-only policy evaluation and an explicit greedy update
for discounted continuous-time LQR.
For this LQR problem, additive noise changes the optimal value only by a
state-independent constant and leaves the optimal feedback gain
unchanged~\citep[Theorem~2]{BijlSchon2019DiscountedLQG}.
Following the deterministic formulation of stochastic LQR
in~\citep[Section~4.3]{ZhuZhangZhang2025OptimalPhiBE}, we therefore consider
the system obtained by setting $\Gamma=0$ in
\eqref{eq:md-lqr-dynamics}, with the same reward and discount.
The batch is still generated by the noisy
system~\eqref{eq:md-lqr-dynamics}; setting $\Gamma=0$ changes only the
evaluation target, not the data.

For a discounted-stabilizing policy $u=-Kx$, let $V_0^K$ denote its value
for this deterministic system. Define the drift-only state--action
critic by $\beta Q_0^K(x,u)=\mathcal R(x,u)+(Ax+Bu)^\top\nabla V_0^K(x)$.
The policy-evaluation equation gives $Q_0^K(x,-Kx)=V_0^K(x)$, so the value
can be recovered from the critic. Policy improvement then maximizes
$Q_0^K(x,u)$ over $u$ at each state.

However, the drift coefficients $A,B$ are unknown. PhiBE replaces the
drift in the evaluation equation by the conditional mean increment
$b_{\dt}(x,u):=\E[(X'-x)/\dt\mid X=x,u]$~\citep{Zhu2024PhiBE}.
For the linear system considered here,
$b_{\dt}(x,u)=\bar A_\Delta x+\bar B_\Delta u$, where
$\bar A_\Delta=(e^{A\Delta t}-I)/\Delta t$ and
$\bar B_\Delta=\Delta t^{-1}\int_0^{\Delta t}e^{A\tau}B\,d\tau$.
These matrices characterize the population approximation; the
sample-based evaluation below uses observed increments directly.

The resulting evaluation equation corresponds to a deterministic
surrogate LQR problem with dynamics $\dot x=b_{\dt}(x,u)$ and the same
reward and discount. For a policy $u=-Kx$ that is discounted-stabilizing
for this surrogate, its value and critic satisfy
\begin{equation}
\label{eq:md-new-Q}
\beta Q_{\dt}^{K}(x,u)=\mathcal R(x,u)
+(\bar A_\Delta x+\bar B_\Delta u)^\top\nabla V_{\dt}^{K}(x), \qquad V_{\dt}^{K}(x)=Q_{\dt}^{K}(x,-Kx).
\end{equation}

The LQR structure makes both evaluation and improvement explicit.
Writing $F_{\Delta,K}:=\bar A_\Delta-\bar B_\Delta K-\beta I/2$,
classical policy evaluation~\citep{Kleinman1968} gives
$V_{\dt}^{K}(x)=-x^\top\Pi_Kx$, where $\Pi_K\succeq0$ is the unique
symmetric solution of
\begin{equation}
\label{eq:md-policy-lyapunov}
F_{\Delta,K}^\top\Pi_K+\Pi_KF_{\Delta,K}+Q+K^\top RK=0.
\end{equation}
Consequently, $Q_{\dt}^{K}$ is quadratic in $(x,u)$, so it can be
represented exactly by quadratic features. It is also strictly concave
in $u$, with Hessian $-2R/\beta$, and its maximizer is the linear feedback
$u=-\mathcal G_\Delta(K)x$, where
$\mathcal G_\Delta(K)=R^{-1}\bar B_\Delta^\top\Pi_K$.
This Hessian does not depend on $\Delta t$: the action
dependence of the critic does not weaken as the observation interval
shrinks.
The full critic expansion is given in \eqref{eq:md-Q-explicit} in
Appendix~\ref{app:quadratic-identities}.

\subsection{One-batch policy iteration}
\label{subsec:lqr-pi-algo}

We implement these evaluation and improvement steps using the same
batch $\mathcal D$ at every iteration.

\paragraph{Policy evaluation.}
We represent the critic as $Q_\theta(x,u)=\Phi(x,u)^\top\theta$, using
the quadratic features
$\Phi(x,u)=(1,\operatorname{svec}(zz^\top)^\top)^\top$, where
$z=(x^\top,u^\top)^\top$. Here $\operatorname{svec}$ stacks the
upper-triangular entries of a symmetric matrix, scaling off-diagonal
entries by $\sqrt2$, and the feature dimension is
$p=1+(d+m)(d+m+1)/2$.
Write $\theta=(a,\operatorname{svec}(H)^\top)^\top$, with $H$ symmetric,
so that $Q_\theta(x,u)=a+z^\top Hz$.

For a fixed gain $K$, the value representation is
$V_\theta^K(x)=Q_\theta(x,-Kx)=\Phi(x,-Kx)^\top\theta$.
For each transition, put $\Phi_i:=\Phi(X_i,u_i)$ and
$v_i:=(X_i'-X_i)/\dt$. Replacing the conditional mean increment in
\eqref{eq:md-new-Q} by $v_i$ gives the sample residual
\[
\delta_{K,i}(\theta)=g_{K,i}^\top\theta-R_i,\qquad g_{K,i}:=\beta\Phi_i-
\left.\frac{\partial}{\partial x}\Phi(x,-Kx)\right|_{x=X_i}v_i.
\]
The derivative acts through both arguments of
$\Phi(x,-Kx)$, with $K$ held fixed.

At the exact surrogate critic, the residual has zero
conditional mean. Following least-squares temporal-difference
learning~\citep{BradtkeBarto1996LSTD}, we enforce
$n^{-1}\sum_i\Phi_i\delta_{K,i}(\theta)=0$ to obtain\footnote{In practice, ridge
regression can be used to regularize policy evaluation.}
\begin{equation}
\label{eq:md-new-residual}
\widehat A_K=\frac1n\sum_i\Phi_i g_{K,i}^\top,\quad
\widehat c=\frac1n\sum_i\Phi_iR_i,\quad
\widehat A_K\widehat\theta_K=\widehat c.
\end{equation}

\paragraph{Policy improvement.}
Let $\widehat H_{K,uu}$ and $\widehat H_{K,ux}$
be the action--action and action--state blocks of the fitted quadratic
matrix $\widehat H_K$.
If $\widehat H_{K,uu}\prec0$\footnote{Under the conditions of Theorem~\ref{thm:large-beta-operator},
part~(i) guarantees that $\widehat A_{K_k}$ is invertible and
$\widehat H_{k,uu}\prec0$ at every iteration.}, the fitted critic has a unique maximizer,
giving the updated gain
$K^+=\widehat H_{K,uu}^{-1}\widehat H_{K,ux}$.

\begin{algorithm2e}[H]
\caption{One-batch PhiBE-PI}\label{alg:one-batch}
\KwIn{$\mathcal D,\beta,\Delta t,K_0$, and iteration count $N$.}
Compute $\widehat c$ using \eqref{eq:md-new-residual}\;
\For{$k=0,\ldots,N-1$}{
Form $\widehat A_{K_k}$ using \eqref{eq:md-new-residual}\;
Solve $\widehat A_{K_k}\widehat\theta_k=\widehat c$ and extract $\widehat H_k$\;
Set $K_{k+1}=\widehat H_{k,uu}^{-1}\widehat H_{k,ux}$\;
}
\KwOut{$K_N$.}
\end{algorithm2e}

\section{Convergence and Accuracy of Policy Iteration}
\label{sec:lqr-pi}

Algorithm~\ref{alg:one-batch} reuses one batch for every evaluation and
greedy update. We first treat the batch as fixed and study how a small
perturbation of the evaluation operator affects convergence and gain
accuracy. Section~\ref{sec:statistical} then bounds this perturbation
with high probability.

Throughout the analysis, the system matrices $A,B$,
cost matrices $Q\succeq0,R\succ0$, and initial gain $K_0$ are fixed.
Vector and coefficient norms are Euclidean; gain and value-matrix
errors use the Frobenius norm.

Let $K^\star_\Delta$ denote the stabilizing optimal gain of the surrogate
LQR problem in Section~\ref{subsec:lqr-q-approximation}.

\subsection{Fixed-operator representation of one-batch PI}
\label{subsec:empirical-quadratic}

When the empirical Gram matrix
$\widehat G=n^{-1}\sum_i\Phi_i\Phi_i^\top$ is invertible, normalizing
the evaluation system~\eqref{eq:md-new-residual} separates its dependence
on the batch from its dependence on the policy.

Let $w$ be the reward coefficient vector, so that
$\mathcal R(x,u)=\Phi(x,u)^\top w$.
For $Q_\theta(x,u)=a+z^\top Hz$, define
$L_K\theta=\operatorname{svec}([I;-K]^\top H[I;-K])$.
Thus, $L_K\theta$ records the quadratic coefficients of
$Q_\theta(x,-Kx)$. The empirical evaluation system satisfies
$\widehat A_K=\widehat G(\beta I-\widehat ML_K),\qquad
\widehat c=\widehat Gw,$
where
\[
\widehat M=\widehat G^{-1}\left[
\frac1n\sum_i\Phi_i
\operatorname{svec}(X_iv_i^\top+v_iX_i^\top)^\top\right].
\]
Hence every evaluation in Algorithm~\ref{alg:one-batch} solves
\begin{equation}
\label{eq:normalized-batch}
(\beta I-\widehat ML_K)\,\widehat\theta_K=w.
\end{equation}
The batch enters only through $\widehat M$, which is fixed across
iterations; the policy enters only through $L_K$.

Whenever its updates are well defined,
Algorithm~\ref{alg:one-batch} performs exact policy iteration for
one fixed empirical operator. Under the conditions of
Theorem~\ref{thm:large-beta-operator}, a single bound on
$\|\widehat M-M_\Delta\|_{\mathrm{op}}$ controls all data-dependent
iterates: estimation error shifts the limiting fixed point without
accumulating across updates. Since $\widehat M$ need not arise from
linear dynamics, classical Riccati theory does not apply directly,
but the quadratic error identity in
Section~\ref{subsec:newton-comparison} survives.

The same construction applies to the population operators $M_\Delta$
and $M_0$ defined in~\eqref{eq:normalized-population}, built from the
drifts $[\bar A_\Delta\ \bar B_\Delta]$ and $[A\ B]$, respectively.
Under the conditions of Theorem~\ref{thm:large-beta-operator}, the
population optima and the empirical limit correspond to the three
operators as follows:
\[
M_0\mapsto K^\star_{\mathrm{ct}},\qquad
M_\Delta\mapsto K^\star_\Delta,\qquad
\widehat M\mapsto\widehat K_\Delta.
\]
Finite-sample error corresponds to the perturbation $\widehat M-M_\Delta$,
while discretization error corresponds to $M_\Delta-M_0$.

\subsection{Quadratic error identities}
\label{subsec:newton-comparison}

With exact population evaluation, $\mathcal G_\Delta$ is Kleinman's
iteration for the surrogate Riccati equation, whose value-matrix
iteration is Newton's method~\citep{Kleinman1968}.
Let $\Pi^\star_\Delta$ be the stabilizing surrogate Riccati solution and
define $\mathcal L_{\Delta,K}(U)=-F_{\Delta,K}^\top U-UF_{\Delta,K}$.
Completing the square gives
\begin{equation}
\label{eq:md-quadratic-identity}
\mathcal L_{\Delta,K}(\Pi_K-\Pi^\star_\Delta)
=(K-K^\star_\Delta)^\top R(K-K^\star_\Delta),
\end{equation}
so the update error $\mathcal G_\Delta(K)-K^\star_\Delta
=R^{-1}\bar B_\Delta^\top(\Pi_K-\Pi^\star_\Delta)$ is quadratic in
$K-K^\star_\Delta$.

The empirical iteration has an analogous identity. Let
$\widehat K_\Delta$ be a fixed point with well-defined evaluation and
greedy improvement, and let $\theta^\dagger$ be its evaluated critic,
with quadratic matrix $H^\dagger$.
The greedy condition $H^\dagger_{ux}=H^\dagger_{uu}\widehat K_\Delta$
cancels the terms linear in $E:=K-\widehat K_\Delta$, and subtracting
the two evaluation equations gives
\begin{equation}
\label{eq:main-empirical-quadratic}
\widehat\theta_K-\theta^\dagger
=(\beta I-\widehat ML_K)^{-1}\widehat M
\operatorname{svec}(E^\top H^\dagger_{uu}E)
\end{equation}
whenever evaluation at $K$ is invertible
(Lemma~\ref{lem:empirical-quadratic},
Appendix~\ref{app:empirical-quadratic}). This identity uses only the
quadratic critic and the greedy condition, not a linear-drift structure
of $\widehat M$, so the quadratic error structure survives the empirical
perturbation.

\subsection{Convergence and gain accuracy}
\label{subsec:fixed-batch-convergence}

The empirical fixed point separates the gain error into three terms:
\begin{equation}
\label{eq:md-three-term-bound}
\|K_N-K^\star_{\mathrm{ct}}\|_F\le
\underbrace{\|K_N-\widehat K_\Delta\|_F}_{\text{iteration}}
+\underbrace{\|\widehat K_\Delta-K^\star_\Delta\|_F}_{\text{finite-sample}}
+\underbrace{\|K^\star_\Delta-K^\star_{\mathrm{ct}}\|_F}_{\text{discretization}}.
\end{equation}

\begin{theorem}[Fixed-batch policy iteration: convergence and gain accuracy]
\label{thm:large-beta-operator}
There exist thresholds $\beta_0,h_0,\varepsilon_0>0$ and constants
$C_{\mathrm{emp}},C_{\mathrm{PI}},C_\star,C_{\mathrm{disc}}>0$,
all depending only on $A,B,Q,R,K_0$, with the
following property. For $\beta\ge\beta_0$ and $0<\Delta t\le h_0$,
suppose the empirical Gram matrix $\widehat G$ is invertible and
the operator estimation error satisfies
$\varepsilon_{\mathrm{op}}:=\|\widehat M-M_\Delta\|_{\mathrm{op}}\le\varepsilon_0$.
Then Algorithm~\ref{alg:one-batch}, starting at $K_0$, satisfies:
\begin{enumerate}[label=(\roman*),leftmargin=*,nosep]
\item \emph{Stability and well-posedness.} Every policy evaluation is
invertible, every evaluated critic has $\widehat H_{K_k,uu}\prec0$,
and every iterate $K_k$ is discounted-stabilizing for both the true
and surrogate systems.
\item \emph{Quadratic convergence.} The iterates converge to an
empirical fixed point $\widehat K_\Delta$. For
$d_k=\|K_k-\widehat K_\Delta\|_F$,
\begin{equation}
\label{eq:md-large-beta-empirical}
d_{k+1}\le\frac{C_{\mathrm{emp}}}{\beta}d_k^2\qquad(k\ge0).
\end{equation}
\item \emph{Gain accuracy.} For every $N\ge1$,
\begin{equation}
\label{eq:md-deterministic-end-to-end}
\|K_N-K^\star_{\mathrm{ct}}\|_F
\le\frac1\beta\left(
C_{\mathrm{PI}}\,2^{-2^{N-1}}
+C_\star\varepsilon_{\mathrm{op}}+C_{\mathrm{disc}}\Delta t \right).
\end{equation}
\end{enumerate}
\end{theorem}

The three terms in
\eqref{eq:md-deterministic-end-to-end} bound, respectively, the iteration
error relative to the empirical fixed point, the finite-sample error
between the empirical and population fixed points, and the discretization
error between the surrogate and continuous-time optima, as decomposed in
\eqref{eq:md-three-term-bound}. For large $\beta$, the normalized
evaluation operator $I-\beta^{-1}\widehat ML_K$ is uniformly close
to the identity and $\beta\widehat H_{K,uu}$ is close to $-R$.
These estimates ensure well-posed evaluation and a contractive greedy map.
The updated gains are $O(1/\beta)$, and an operator perturbation of
size $\delta$ shifts the fixed point by $O(\delta/\beta)$, accounting
for the common factor $1/\beta$ in
\eqref{eq:md-deterministic-end-to-end}.

The quadratic recursion
\eqref{eq:md-large-beta-empirical} yields doubly exponential decay of
the iteration error. Once the iteration term is small, the finite-sample
and discretization terms dominate the accuracy bound.
Section~\ref{sec:statistical} bounds $\varepsilon_{\mathrm{op}}$ in
terms of sample size and training interval.
Appendix~\ref{app:convergence-experiment} illustrates convergence to the
empirical fixed point, with observed local orders approximately two.
The proof of Theorem~\ref{thm:large-beta-operator} is given in
Appendix~\ref{app:large-beta}.

\section{Finite-Sample Control Guarantees}
\label{sec:statistical}

Theorem~\ref{thm:large-beta-operator} reduces the
analysis to two random quantities: the empirical Gram matrix
$\widehat G$ and the operator estimation error $\varepsilon_{\mathrm{op}}$.
Exact quadratic closure gives
\begin{equation}
\label{eq:normalized-cancellation}
\widehat M-M_\Delta=\widehat G^{-1}\widehat W,
\end{equation}
where $\widehat W$ is an empirical noise matrix formed
from the deviations of the normalized increments from their conditional
means; see \eqref{eq:md-lift-noise}. Neither $\widehat G$ nor
$\widehat W$ depends on the policy, so concentrating these two matrices
gives one event that covers every iterate, including policies selected
using the batch. The normalized increment noise has standard deviation
$O(\Delta t^{-1/2})$. Under the sampling and coverage conditions below,
this leads to an operator-error bound of order $(n\Delta t)^{-1/2}$,
up to confidence factors.

\begin{assumption}[Dynamics, initialization, and data]
\label{ass:global-data}
\begin{enumerate}[label=(\roman*),leftmargin=*]
\item \emph{Dynamics and initialization.}
The pair $(A-\beta I/2,B)$ is stabilizable and
$(Q^{1/2},A-\beta I/2)$ is detectable. The fixed initial gain $K_0$ is
independent of the batch and $A-BK_0-\beta I/2$ is Hurwitz.
\item \emph{Sampling.}
Starting pairs $(X_i,u_i)$ are i.i.d. from a fixed distribution $\mu$,
with $\|X_i\|\le M_x$ and $\|\Phi_i\|\le M_\Phi$.
Successors follow the held-action data model in
Section~\ref{subsec:lqr-1d-setup}. The transition noises are independent
across samples and of the starting pairs.
\item \emph{Coverage.}
The population Gram matrix satisfies
$G_\Phi:=\E[\Phi_i\Phi_i^\top]\succeq\lambda_\Phi I\succ0$.
\end{enumerate}
\end{assumption}

Coverage is a rank condition on quadratic features; it does not require
samples near the optimal policy. No positive-definite diffusion covariance,
full actuation, or positive-definite $Q$ is required.

\begin{theorem}[Finite-sample gain accuracy]
\label{cor:md-large-beta-end-to-end}
Suppose Assumption~\ref{ass:global-data} holds.
There exist positive constants $\beta_0,h_0,c_1,c_2$ such that, for every
$\beta\ge\beta_0$, $0<\Delta t\le h_0$, and $\delta\in(0,1)$, if
\begin{equation}
\label{eq:md-sample-size-conditions}
n\ge c_1\log(4p/\delta),\qquad
n\Delta t\ge c_2\log(4p/\delta),
\end{equation}
then, with probability at least $1-\delta$, the conditions of
Theorem~\ref{thm:large-beta-operator} hold. On this same event,
simultaneously for every $N\ge1$,
\begin{equation}
\label{eq:md-large-beta-end-to-end}
\|K_N-K^\star_{\mathrm{ct}}\|_F
\le\frac1\beta\left(C_{\mathrm{PI}}\,2^{-2^{N-1}}
+C_{\mathrm{stat}}\sqrt{\frac{\log(4p/\delta)}{n\Delta t}}
+C_{\mathrm{disc}}\Delta t\right).
\end{equation}
Here $C_{\mathrm{PI}},C_{\mathrm{stat}},C_{\mathrm{disc}}>0$.
All constants depend only on the fixed system,
cost, initializer, sampling bounds, and coverage constant, and are
independent of $\beta,\Delta t,n,N,\delta$.
\end{theorem}

The proof of
Theorem~\ref{cor:md-large-beta-end-to-end} is given in
Appendix~\ref{app:finite-sample-gain-proof}.

\begin{remark}[Observation interval and sample size]
In \eqref{eq:md-sample-size-conditions}, $n$
controls Gram invertibility and $n\Delta t$ controls the operator error.
At fixed $n$, reducing $\Delta t$ lowers discretization bias but amplifies
noise in the normalized increments. For $\beta\ge\beta_0$, balancing
the finite-sample and discretization terms gives $\Delta t\asymp n^{-1/3}$,
independently of $\beta$ because both terms share $1/\beta$.
Once the iteration term is negligible, the bound gives gain error
$O_p(\beta^{-1}n^{-1/3})$ and, for a fixed problem, cost gap
$O_p(n^{-2/3})$ (Corollary~\ref{cor:md-cost}).
Table~\ref{tab:scaling-slopes} in Appendix~\ref{app:scaling-details}
reports slopes close to $-1/3$ and $-2/3$, respectively.
\end{remark}

\begin{remark}[Small and moderate discounts]
\label{rem:small-beta}
For any fixed $\beta>0$ under
Assumption~\ref{ass:global-data}, global convergence of Kleinman's
iteration and a perturbation argument yield, with high probability,
discounted stability of all iterates, quadratic convergence after finitely
many steps, and a gain bound with the same three terms as
\eqref{eq:md-large-beta-end-to-end}. We omit this bound because its
constants depend implicitly on $\beta$ and $K_0$, obscuring the discount
dependence. The explicit $1/\beta$ scaling should not be extrapolated
to small $\beta$: Figure~\ref{fig:beta-dt-control}(a) no longer follows
it below $\beta\approx1$, and the scalar discretization error is
proportional to $\beta\Delta t$ for small $\beta$
\citep[cf.][Theorem~4.4]{ZhuZhangZhang2025OptimalPhiBE}.
\end{remark}

\section{Numerical Experiments}
\label{subsec:lqr-numerics-main}

\begin{figure}[!htbp]
\centering
\includegraphics[width=\linewidth]{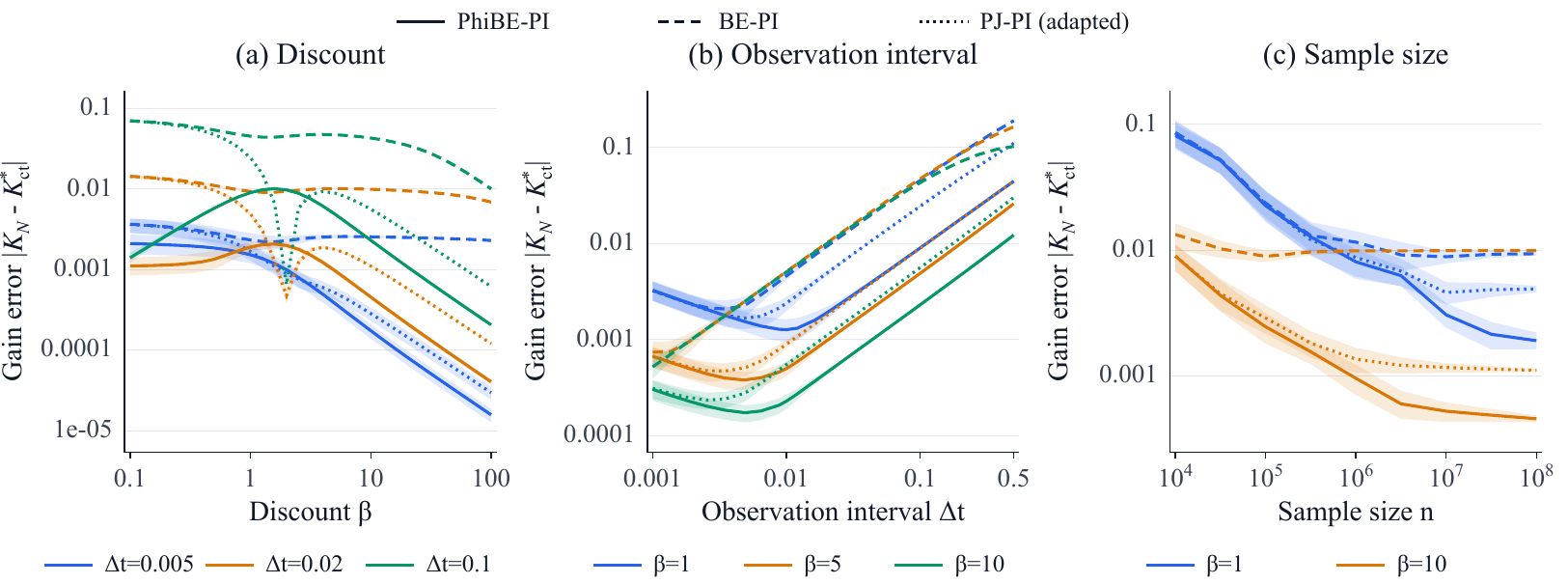}
\caption{Gain accuracy in the scalar system with process
noise. (a,b) Gain error
$|K_N-K^\star_{\mathrm{ct}}|$ at $n=10^8$;
colors give the fixed $\Delta t$ or $\beta$ listed below each panel.
(c) Gain error versus sample size at fixed $\Delta t=0.02$.
Curves and bands show means and pointwise 95\% bootstrap
confidence intervals over 30 paired batches.}
\label{fig:beta-dt-control}
\end{figure}

We compare PhiBE-PI with two discrete-time Bellman-equation baselines
that reuse the same batch and initial gain: BE-PI, with discount
$e^{-\beta\Delta t}$ and stage reward $\Delta t\,\mathcal R(x,u)$,
and an adaptation of \citet{PangJiang2023Stationary} to endpoint
transitions, denoted PJ-PI (adapted), which is equivalent to BE-PI with
discount $(1+\beta\Delta t)^{-1}$ (Appendix~\ref{app:pj-adaptation}).
This baseline also corresponds to the modified
Optimal-BE formulation in
\citet[Eq.~(38)]{ZhuZhangZhang2025OptimalPhiBE}.
The scalar system has $A=0.4$,
$B=\Gamma=Q=R=1$, and $K_0=1.5$; Appendix~\ref{app:additional-experiments}
gives data-generation protocols and noisy multidimensional runs.

\paragraph{Discount, observation interval, and sample size.}
Figure~\ref{fig:beta-dt-control} examines the three
parameter dependencies in Theorem~\ref{cor:md-large-beta-end-to-end}.
Panel~(a) shows that, for large discounts, the gain error of PhiBE-PI
decreases approximately as $1/\beta$. BE-PI has substantially larger
errors when $\beta\Delta t$ is not small, while PJ-PI (adapted) also
maintains small errors.

Panel~(b) shows the predicted tradeoff: reducing
$\Delta t$ lowers discretization bias but amplifies noise in normalized
state increments, producing U-shaped gain-error curves.

Panel~(c) shows how increasing the sample size reduces
finite-sample error until the total error reaches a plateau determined
by each method's population bias. At $\beta=10$, PhiBE-PI reaches the
lowest plateau among the three methods.

\begin{figure}[!htbp]
\centering
\includegraphics[width=\linewidth]{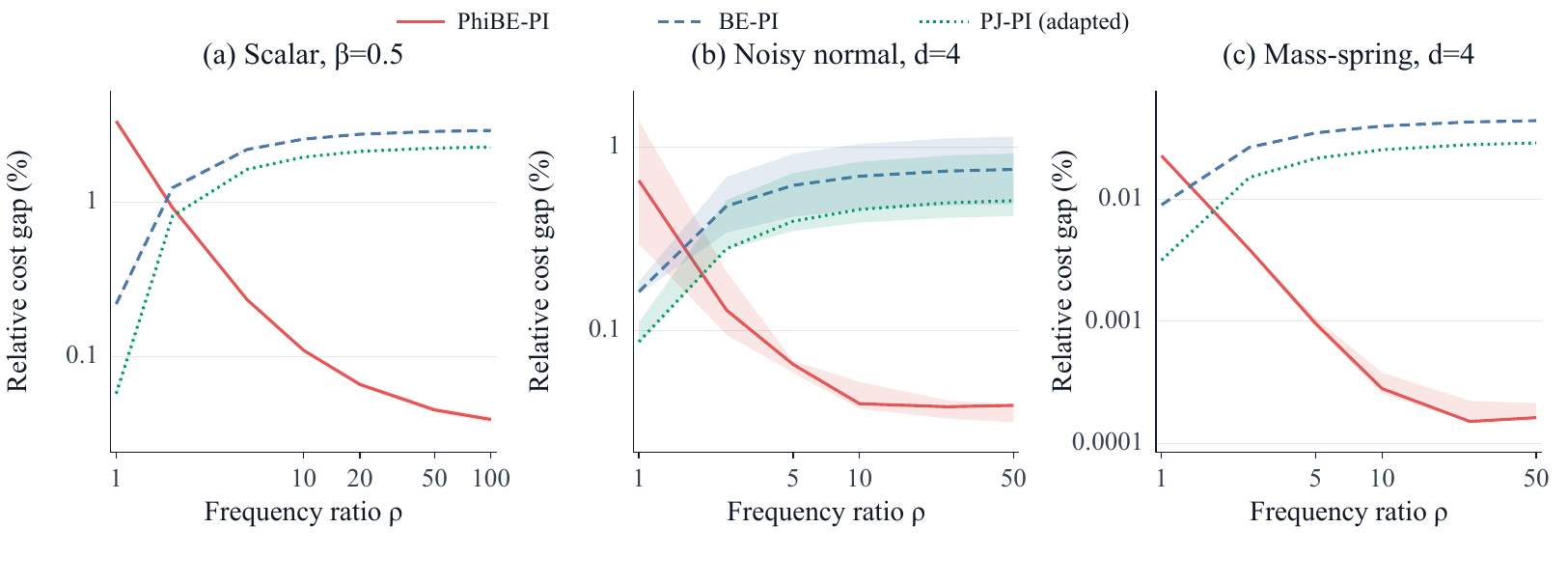}
\caption{Relative cost gaps for gains trained at
$\Delta_{\mathrm{tr}}=0.5$ and deployed at
$\Delta_{\mathrm{ts}}=\Delta_{\mathrm{tr}}/\rho$, using the optimal
held-action policy at each deployment interval as the reference.
(a) Noisy scalar system: means and 95\% bootstrap
confidence intervals. (b) Noisy normal systems and (c) overdamped
mass--spring systems: medians and interquartile bands across five systems,
after averaging batches within each system. Protocols are in
Appendix~\ref{app:multid-deployment}.}
\label{fig:deployment-ratio-control}
\end{figure}

\paragraph{Training and deployment rates.}
Figure~\ref{fig:deployment-ratio-control} fixes the
gains learned at $\Delta_{\mathrm{tr}}=0.5$ and varies the deployment interval
$\Delta_{\mathrm{ts}}=\Delta_{\mathrm{tr}}/\rho$.
In the scalar system, the baselines perform better at $\rho=1$, but
PhiBE-PI becomes more accurate as deployment becomes faster. At
$\rho=100$, its relative cost gap is $0.04\%$, compared with
$2.92\%$ for BE-PI and $2.28\%$ for PJ-PI (adapted). This reversal
reflects their different targets: the baselines optimize sampled MDPs
tied to the training interval, whereas PhiBE-PI approximates the
continuous-time optimum.

Panels~(b,c) show similar benefits on noisy four-state
normal systems and overdamped mass--spring chains. The lightly damped
example in Appendix~\ref{app:lightly-damped-mass-spring} illustrates
a limitation: at a coarse training interval, the first-order drift
approximation introduces substantial policy bias, and PhiBE-PI
performs worse than both baselines even with noiseless data.
Faster deployment cannot remove this training-time bias.

\FloatBarrier

\section{Conclusion}
\label{sec:conclusion}

We developed PhiBE-PI for continuous-time LQR from a single batch of
noisy discrete-time off-policy transitions. A shared policy-independent
empirical operator allows one concentration event to control all
data-dependent iterates. Under the stated sampling assumptions,
sufficiently large discounts, small training intervals, and sufficient
data ensure discounted stability and quadratic convergence to an
empirical fixed point. The gain bound separates iteration, finite-sample,
and discretization errors, with the latter two bounded at orders
$1/(\beta\sqrt{n\Delta t})$ and $\Delta t/\beta$. Experiments illustrate
the tradeoff in choosing the training interval and the benefit of faster
deployment when the training interval is small relative to the system's
time scale.

Future work will focus on systematic benchmarking across system dynamics,
discount rates, training intervals, and sample budgets to identify when
PhiBE-PI improves on discrete-time baselines and which error mechanisms
limit its performance.

\subsection*{Generative AI Use}
Generative AI tools assisted with literature review, manuscript drafting and revision, proof development, numerical checks, and LaTeX preparation.

\bibliography{\preprintbib}
\clearpage
\appendix
\section{Quadratic Evaluation and Operator Identities}
\label{app:quadratic-identities}

Substituting $V_{\dt}^K(x)=-x^\top\Pi_Kx$ into \eqref{eq:md-new-Q} gives the quadratic critic
\begin{equation}
\label{eq:md-Q-explicit}
\begin{split}
Q_{\dt}^{K}(x,u)=-\frac1\beta\bigl[&x^\top(Q+\bar A_\Delta^\top\Pi_K+
\Pi_K\bar A_\Delta)x+2x^\top\Pi_K\bar B_\Delta u+u^\top Ru\bigr].
\end{split}
\end{equation}

\subsection{Quadratic evaluation and closure}

\paragraph{Derivation and uniqueness of \eqref{eq:md-Q-explicit}.}
For any quadratic critic, write its restriction to $u=-Kx$ as
$V_{\dt}^{K}(x)=c+\ell^\top x-x^\top\Pi x$, with $\Pi$ symmetric.
Restricting \eqref{eq:md-new-Q} to this policy and matching constant
and linear terms gives $\beta c=0$ and
$(F_{\Delta,K}^\top-\beta I/2)\ell=0$.
Since $\beta>0$ and $F_{\Delta,K}$ is Hurwitz, these imply
$c=0$ and $\ell=0$. Matching quadratic terms gives
\eqref{eq:md-policy-lyapunov}. Hurwitz stability makes its symmetric
solution unique, and $Q+K^\top RK\succeq0$ implies $\Pi=\Pi_K\succeq0$
by the Lyapunov integral representation.
Substitution into \eqref{eq:md-new-Q} now uniquely determines the critic
as \eqref{eq:md-Q-explicit}, with no constant or linear terms.
Finally, its Hessian in $u$ is $-2R/\beta\prec0$, and setting its
gradient to zero yields the unique maximizer
$u=-R^{-1}\bar B_\Delta^\top\Pi_Kx$.

\paragraph{Population quadratic identity.}
Let $K^\star_\Delta=R^{-1}\bar B_\Delta^\top\Pi^\star_\Delta$.
Rewriting the surrogate CARE with $F_{\Delta,K}$ gives
\[
F_{\Delta,K}^\top\Pi^\star_\Delta+\Pi^\star_\Delta F_{\Delta,K}
+Q+K^\top RK
=(K-K^\star_\Delta)^\top R(K-K^\star_\Delta).
\]
Subtracting this from \eqref{eq:md-policy-lyapunov} gives
\eqref{eq:md-quadratic-identity}.

\subsection{Normalized operators and shared noise}
\label{app:normalized-operator}
Write $S_K=[I;-K]$, so that the quadratic restriction is
$L_K\theta=\operatorname{svec}(S_K^\top HS_K)$.
The directional derivative satisfies
\[
\theta^\top\left[\frac{\partial}{\partial x}\Phi(x,-Kx)\right]v
=(L_K\theta)^\top\operatorname{svec}(xv^\top+vx^\top).
\]
Define
\[
\widehat V=\frac1n\sum_i\Phi_i
\operatorname{svec}(X_iv_i^\top+v_iX_i^\top)^\top.
\]
When $\widehat G$ is invertible, $\widehat M=\widehat G^{-1}\widehat V$
and the empirical evaluation system factors as
$\widehat A_K=\beta\widehat G-\widehat V L_K
=\widehat G(\beta I-\widehat M L_K)$ and $\widehat c=\widehat G w$,
where $w=(0,-\operatorname{svec}(\operatorname{diag}(Q,R))^\top)^\top$.
For any fixed operator $M$, write
\begin{equation}
\label{eq:main-shared-operator}
C_M(K)\theta_M(K)=w,\qquad C_M(K)=\beta I-ML_K.
\end{equation}
Let $H_{M,K}$ be the quadratic matrix of $\theta_M(K)$, and set
$g_M(K)=H_{M,K,uu}^{-1}H_{M,K,ux}$.
We call $K$ feasible when $C_M(K)$ is invertible and $H_{M,K,uu}\prec0$.
For $M=\widehat M$, these are the evaluation and greedy maps of
Algorithm~\ref{alg:one-batch}.

Put $D_\Delta=[\bar A_\Delta\ \bar B_\Delta]$, $J=[I_d\ 0]$, and define
\begin{equation}
\label{eq:normalized-population}
M_\Delta\operatorname{svec}(U)
=\begin{pmatrix}0\\\operatorname{svec}(D_\Delta^\top UJ+J^\top UD_\Delta)\end{pmatrix}.
\end{equation}
Its continuous extension at $\Delta t=0$ is $M_0$, with $D_0=[A\ B]$.
For a surrogate-stabilizing gain, $C_{M_\Delta}(K)$ is invertible:
if $C_{M_\Delta}(K)\theta=0$, restriction to $u=-Kx$ gives $c=0$ and
$F_{\Delta,K}^\top(S_K^\top HS_K)+(S_K^\top HS_K)F_{\Delta,K}=0$.
Hurwitz stability implies $S_K^\top HS_K=0$, and the full coefficient
equation then gives $H=0$. Thus population evaluation recovers
\eqref{eq:md-Q-explicit} exactly.

Define the increment noise
$\zeta_i=v_i-\bar A_\Delta X_i-\bar B_\Delta u_i$ and set
\begin{equation}
\label{eq:md-lift-noise}
Y_i=\operatorname{svec}(X_i\zeta_i^\top+\zeta_iX_i^\top),\qquad
\widehat W=\frac1n\sum_i\Phi_iY_i^\top.
\end{equation}
The observed quadratic increment satisfies
\[
\operatorname{svec}(X_iv_i^\top+v_iX_i^\top)=M_\Delta^\top\Phi_i+Y_i.
\]
Hence $\widehat V=\widehat G M_\Delta+\widehat W$, proving
\eqref{eq:normalized-cancellation} whenever $\widehat G$ is invertible.
The noise matrix
$\widehat W\in\mathbb R^{p\times d(d+1)/2}$ is independent of the policy.

\subsection{Empirical quadratic convergence}
\label{app:empirical-quadratic}
\begin{lemma}[Exact quadratic identity for a fixed operator]
\label{lem:empirical-quadratic}
Let $K^\dagger=g_M(K^\dagger)$ be feasible for $M$, and let $H^\dagger$
be the quadratic matrix of $\theta^\dagger=\theta_M(K^\dagger)$.
For any other feasible $K$, with $E=K-K^\dagger$,
\begin{equation}
\label{eq:empirical-quadratic-evaluation}
\theta_M(K)-\theta^\dagger
=(\beta I-ML_K)^{-1}M\operatorname{svec}(E^\top H^\dagger_{uu}E).
\end{equation}
In particular, the derivative of $g_M$ vanishes at $K^\dagger$, and
the iteration $K\mapsto g_M(K)$ converges locally quadratically to
$K^\dagger$.
\end{lemma}

\begin{proof}
The two evaluated critics satisfy
\[
C_M(K)\theta_M(K)=w,\qquad
C_M(K^\dagger)\theta^\dagger=w.
\]
Subtracting these equations and collecting the critic
difference on the left gives
\[
C_M(K)(\theta_M(K)-\theta^\dagger)
=M(L_K-L_{K^\dagger})\theta^\dagger.
\]
The fixed-point condition is
$H^\dagger_{ux}=H^\dagger_{uu}K^\dagger$, and symmetry gives
$H^\dagger_{xu}=(K^\dagger)^\top H^\dagger_{uu}$.
Writing $K=K^\dagger+E$, the quadratic restriction expands as
\[
\begin{aligned}
&S_K^\top H^\dagger S_K-S_{K^\dagger}^\top H^\dagger S_{K^\dagger}\\
&\quad=-\bigl(H^\dagger_{xu}-(K^\dagger)^\top H^\dagger_{uu}\bigr)E-E^\top\bigl(H^\dagger_{ux}-H^\dagger_{uu}K^\dagger\bigr)
+E^\top H^\dagger_{uu}E\\
&\quad=E^\top H^\dagger_{uu}E.
\end{aligned}
\]
Consequently,
$(L_K-L_{K^\dagger})\theta^\dagger
=\operatorname{svec}(E^\top H^\dagger_{uu}E)$.
Substituting this into the difference equation and multiplying by
$C_M(K)^{-1}$ proves \eqref{eq:empirical-quadratic-evaluation}.
The greedy map is smooth wherever the action
block is invertible, so local inverse and curvature bounds give
$\|g_M(K)-K^\dagger\|_F\le C\|E\|_F^2$ for a locally bounded $C$.
The conversion below gives the explicit multiplier.
\end{proof}

\paragraph{Greedy-map conversion.}
\label{app:greedy-conversion}
In the notation of Lemma~\ref{lem:empirical-quadratic}, subtracting
$H_{M,K,uu}K^\dagger$ in the greedy formula gives
\begin{equation}
\label{eq:empirical-quadratic-gain}
\begin{split}
g_M(K)-K^\dagger
=H_{M,K,uu}^{-1}\bigl[&(H_{M,K,ux}-H^\dagger_{ux})-(H_{M,K,uu}-H^\dagger_{uu})K^\dagger\bigr].
\end{split}
\end{equation}
For example, the product
$\|H_{M,K,uu}^{-1}\|_{\mathrm{op}}(1+\|K^\dagger\|_{\mathrm{op}})\times\|C_M(K)^{-1}\|_{\mathrm{op}}\|M\|_{\mathrm{op}}\|H^\dagger_{uu}\|_{\mathrm{op}}$
is a valid multiplier of $\|E\|_F^2$. All factors are locally bounded.

\section{Proofs of the Control Guarantees}
\label{app:lqr-multidimensional}

\subsection{Large-discount convergence and perturbation bounds}
\label{app:large-beta-results}
\label{app:large-beta}
\begin{proof}[Proof of Theorem~\ref{thm:large-beta-operator}]
All constants in this subsection may depend on the fixed system, cost,
and initializer, but not on $\beta,\Delta t$ or the batch.
The notation $C$ denotes a finite constant whose value may change
between estimates.

\paragraph{Step 1: Uniform feasibility and contraction.}
Choose a closed Frobenius ball
$\mathcal U=\{K:\|K\|_F\le R_K\}$ with
$R_K>\|K_0\|_F+1$, and fix a small $h_0>0$.
The matrices $\bar A_\Delta,\bar B_\Delta,M_\Delta$ extend continuously
to $\Delta t=0$ and are bounded on $[0,h_0]$.
Choose a fixed $m$ such that the convex operator ball
$\mathcal M=\{M:\|M\|_{\mathrm{op}}\le m\}$ contains $M_0$ and every
$M_\Delta$, with room for a fixed positive perturbation radius.
On $\mathcal U$, both $L_K$ and its derivative are bounded; write
$t=\sup_{\mathcal U}\|L_K\|_{\mathrm{op}}$.

Take $\beta_0\ge2mt$. For $\beta\ge\beta_0$,
the Neumann series gives
\begin{equation}
\label{eq:large-beta-inverse}
\|C_M(K)^{-1}\|_{\mathrm{op}}\le\frac2\beta,
\qquad
\|\theta_M(K)\|_2\le\frac{2\|w\|_2}{\beta}
\end{equation}
uniformly on $\mathcal M\times\mathcal U$.
From $\beta\theta_M=w+ML_K\theta_M$ and the reward coefficients,
\[
\beta H_{M,K}=-\operatorname{diag}(Q,R)+O(\beta^{-1}).
\]
Here and below the remainder bounds are uniform on the stated sets.
Increasing $\beta_0$ therefore gives
\begin{equation}
\label{eq:large-beta-curvature}
H_{M,K,uu}\preceq-\frac{R}{2\beta},\qquad
\|H_{M,K,uu}^{-1}\|_{\mathrm{op}}\le
\frac{2\beta}{\lambda_{\min}(R)}.
\end{equation}
Since the reward has no $ux$ block,
$H_{M,K,ux}=O(\beta^{-2})$, and hence
\begin{equation}
\label{eq:large-beta-map-bounds}
\|g_M(K)\|_F\le\frac{C_0}{\beta}.
\end{equation}
Differentiation of the evaluation equation gives
\[
D_K\theta_M(K)[E]
=C_M(K)^{-1}M(D_KL_K[E])\theta_M(K).
\]
Its norm is at most $C\|E\|_F/\beta^2$.
Differentiating $H_{M,K,uu}g_M(K)=H_{M,K,ux}$ and using
\eqref{eq:large-beta-curvature} yields
$\|D_Kg_M(K)\|_{F\to F}\le C/\beta$.
We may thus choose $\beta_0$ so that every $g_M$ maps
$\mathcal U$ into its interior and is a contraction there with
Lipschitz constant at most $1/2$.
It has a unique fixed point $K^\dagger(M)$ in $\mathcal U$, and
$\|K^\dagger(M)\|_F\le C_0/\beta$.

\paragraph{Step 2: Discounted stability and population fixed points.}
The population update $\mathcal G_\Delta$ is the
Kleinman update. It can also be obtained through a particular
quasi-Newton stepsize in model-based continuous-time
LQR~\citep[Section~6.1]{BuMesbahiMesbahi2020ContinuousLQR}.
Let $\operatorname{sym}(D)=(D+D^\top)/2$ and set
\[
a=\sup_{\substack{K\in\mathcal U\\0\le\Delta t\le h_0}}
\lambda_{\max}\!\left(
\operatorname{sym}(\bar A_\Delta-\bar B_\Delta K)\right),
\]
where the endpoint $\Delta t=0$ uses $A,B$.
For $\beta>2a$, every gain in $\mathcal U$ is discounted-stabilizing
for both the original and surrogate systems. Increasing $\beta_0$
if necessary, the Lyapunov integral representation gives
\begin{equation}
\label{eq:large-beta-lyapunov-inverse}
\|\mathcal L_{\Delta,K}^{-1}\|_{F\to F}
\le\int_0^\infty e^{-(\beta-2a)s}\,ds
=\frac1{\beta-2a}\le\frac2\beta.
\end{equation}
In particular,
$\|\Pi_K\|_F\le C/\beta$ uniformly.
At the population fixed point, the policy Lyapunov equation and
$K^\dagger=R^{-1}\bar B_\Delta^\top\Pi_{K^\dagger}$ give the
surrogate CARE. The gain is discounted-stabilizing, so this is the
stabilizing optimal solution. Thus
$K^\dagger(M_\Delta)=K^\star_\Delta$ and
$K^\dagger(M_0)=K^\star_{\mathrm{ct}}$.
Moreover, for these large discounts the shifted open-loop
matrices are themselves Hurwitz, so stabilizability and detectability
hold without additional restrictions on $A,B,Q$.

\paragraph{Step 3: Empirical quadratic convergence.}
Choose $\varepsilon_0>0$ so that
$\|M-M_\Delta\|_{\mathrm{op}}\le\varepsilon_0$ implies $M\in\mathcal M$
for every $0\le\Delta t\le h_0$. The feasibility, invariance, and
contraction estimates above apply to $M=\widehat M$.
The normalized evaluation is exactly Algorithm~\ref{alg:one-batch} when
$\widehat G$ is invertible. Its iterates remain in $\mathcal U$, where
both systems are discounted-stabilizing, every evaluation is invertible,
and every action block is negative definite. They converge to the
unique fixed point $\widehat K_\Delta=K^\dagger(\widehat M)$ in
$\mathcal U$.

Taking $M=\widehat M$ and
$K^\dagger=\widehat K_\Delta$ in
Lemma~\ref{lem:empirical-quadratic} proves
\eqref{eq:main-empirical-quadratic}. The evaluator inverse is
$O(\beta^{-1})$, the action-block inverse is $O(\beta)$,
$\|H^\dagger_{uu}\|_{\mathrm{op}}=O(\beta^{-1})$, and $M$ is bounded.
Its multiplier of $\|K-K^\dagger\|_F^2$ is therefore at most
$C_{\mathrm{emp}}/\beta$, uniformly for $K\in\mathcal U$.
This proves \eqref{eq:md-large-beta-empirical}.
The initial error $d_0$ is uniformly bounded, so
$d_1\le D_1/\beta$ for a fixed $D_1>0$.
Increase $\beta_0$ if necessary so that
$C_{\mathrm{emp}}D_1/\beta_0^2\le1/2$.
Repeated squaring gives, for $N\ge1$,
\[
\begin{aligned}
d_N
&\le d_1
\left(\frac{C_{\mathrm{emp}}d_1}{\beta}\right)^{2^{N-1}-1}
\le\frac{2D_1}{\beta}\,2^{-2^{N-1}}.
\end{aligned}
\]
The case $d_1=0$ is immediate. Taking $C_{\mathrm{PI}}=2D_1$
gives
$d_N\le (C_{\mathrm{PI}}/\beta)\,2^{-2^{N-1}}$.

\paragraph{Step 4: Finite-sample and discretization perturbation bounds.}
For $M_1,M_2\in\mathcal M$, subtracting their evaluation equations gives
\[
\theta_{M_1}(K)-\theta_{M_2}(K)
=C_{M_1}(K)^{-1}(M_1-M_2)L_K\theta_{M_2}(K).
\]
Thus the coefficient difference is bounded by
$C\|M_1-M_2\|_{\mathrm{op}}/\beta^2$.
Subtracting the two greedy equations then gives
\begin{equation}
\label{eq:large-beta-map-sensitivity}
\sup_{K\in\mathcal U}\|g_{M_1}(K)-g_{M_2}(K)\|_F
\le\frac{C}{\beta}\|M_1-M_2\|_{\mathrm{op}}.
\end{equation}
The contraction bound implies the corresponding fixed-point estimate
\begin{equation}
\label{eq:large-beta-root-sensitivity}
\|K^\dagger(M_1)-K^\dagger(M_2)\|_F
\le\frac{C_\star}{\beta}\|M_1-M_2\|_{\mathrm{op}}.
\end{equation}

Applying \eqref{eq:large-beta-root-sensitivity} with
$M_1=\widehat M$ and $M_2=M_\Delta$ gives
\begin{equation}
\label{eq:md-large-beta-greedy-error}
\|\widehat K_\Delta-K^\star_\Delta\|_F
\le\frac{C_\star}{\beta}\varepsilon_{\mathrm{op}}.
\end{equation}

Finally, the smooth surrogate coefficients satisfy
$\|M_\Delta-M_0\|_{\mathrm{op}}\le C\Delta t$ on $[0,h_0]$.
Apply \eqref{eq:large-beta-root-sensitivity} to the two population
operators to obtain
\begin{equation}
\label{eq:md-large-beta-components}
\|K^\star_\Delta-K^\star_{\mathrm{ct}}\|_F
\le\frac{C_{\mathrm{disc}}}{\beta}\Delta t.
\end{equation}
Combining the iteration bound with
\eqref{eq:md-large-beta-greedy-error}
and~\eqref{eq:md-large-beta-components} through
\eqref{eq:md-three-term-bound} proves
\eqref{eq:md-deterministic-end-to-end}.
All thresholds were chosen from fixed bounded sets and remain valid
for every $\beta\ge\beta_0$.
\end{proof}

\subsection{Concentration of the empirical operator}
Under the held input, the observed increments satisfy
\begin{equation}
\label{eq:md-transition-noise}
v_i=\bar A_\Delta X_i+\bar B_\Delta u_i+\zeta_i,\quad
\zeta_i\sim\mathcal N(0,\bar\Sigma_\Delta/\Delta t),
\end{equation}
where $\bar\Sigma_\Delta=(\Delta t)^{-1}\int_0^{\Delta t}
 e^{A\tau}\Gamma\Gamma^\top e^{A^\top\tau}\,d\tau$.
The noises are independent across transitions and of the starting pairs.

For the shared noise matrix $\widehat W$ in
\eqref{eq:md-lift-noise}, let $\tau_\Delta=\sqrt{\operatorname{tr}\bar\Sigma_\Delta}$ and
$\sigma_\Delta=\sqrt{\|\bar\Sigma_\Delta\|_{\mathrm{op}}}$.
For $\delta\in(0,1)$, define the concentration radii
\begin{equation}
\label{eq:md-concentration-radii}
\begin{split}
r_G&=M_\Phi^2\left[
\sqrt{\frac{2\log(4p/\delta)}n}+\frac{2\log(4p/\delta)}{3n}\right],\\
r_W&=\frac{2M_\Phi M_x}{\sqrt{n\Delta t}}
\left[\tau_\Delta+\sigma_\Delta\sqrt{2\log(2/\delta)}\right].
\end{split}
\end{equation}
The corresponding operator-error bound is
$\varepsilon_{\mathrm{op},n,\Delta}:=2r_W/\lambda_\Phi$.
Since $\tau_\Delta\le\sqrt d\,\sigma_\Delta$, the bound
depends on the dimension through $\tau_\Delta$.

\begin{proposition}[Empirical operator concentration]
\label{prop:operator-concentration}
Under the sampling and coverage conditions of
Assumption~\ref{ass:global-data}, with probability at least $1-\delta$,
\[
\|\widehat G-G_\Phi\|_{\mathrm{op}}\le r_G,\qquad
\|\widehat W\|_{\mathrm{op}}\le\|\widehat W\|_F\le r_W.
\]
If $r_G\le\lambda_\Phi/2$, then on the same event $\widehat G$ is
invertible and $\|\widehat M-M_\Delta\|_{\mathrm{op}}
\le2r_W/\lambda_\Phi=\varepsilon_{\mathrm{op},n,\Delta}$.
\end{proposition}

\begin{proof}
For $Z_i=\Phi_i\Phi_i^\top-G_\Phi$,
$\|Z_i\|_{\mathrm{op}}\le M_\Phi^2$ and
$\|\E Z_i^2\|_{\mathrm{op}}\le M_\Phi^4$.
Matrix Bernstein gives the stated Gram bound with failure probability
at most $\delta/2$.
Conditional on all state--action pairs, $\operatorname{vec}\widehat W$
is a centered Gaussian vector. Using the isometry of $\operatorname{svec}$,
\[
\|x\zeta^\top+\zeta x^\top\|_F\le2\|x\|\|\zeta\|.
\]
Thus the trace and largest eigenvalue of its conditional covariance are
bounded, respectively, by
\[
\frac{4M_\Phi^2M_x^2\tau_\Delta^2}{n\Delta t},\qquad
\frac{4M_\Phi^2M_x^2\sigma_\Delta^2}{n\Delta t}.
\]
For a centered Gaussian vector with covariance $V$, its norm is at most
$\sqrt{\operatorname{tr}V}+
\sqrt{2\|V\|_{\mathrm{op}}\log(2/\delta)}$ except on an event of
probability $\delta/2$. This proves the bound $r_W$, first conditionally
and then unconditionally. A union bound gives the common event.

If $r_G\le\lambda_\Phi/2$, this event also gives
$\lambda_{\min}(\widehat G)\ge\lambda_\Phi/2$, so
\eqref{eq:normalized-cancellation} implies
\[
\|\widehat M-M_\Delta\|_{\mathrm{op}}
\le\|\widehat G^{-1}\|_{\mathrm{op}}\|\widehat W\|_{\mathrm{op}}
\le\frac{2r_W}{\lambda_\Phi}=\varepsilon_{\mathrm{op},n,\Delta}.
\]

\end{proof}

\subsection{Proof of the finite-sample gain theorem}
\label{app:finite-sample-gain-proof}
\begin{proof}[Proof of Theorem~\ref{cor:md-large-beta-end-to-end}]
Take $\beta_0,h_0,\varepsilon_0,C_{\mathrm{PI}},C_\star,
C_{\mathrm{disc}}$ from
Theorem~\ref{thm:large-beta-operator}, and put
$L_\delta=\log(4p/\delta)>1$.
For $0<\Delta t\le h_0$, the covariance formula in
\eqref{eq:md-transition-noise} gives
\[
\tau_\Delta\le\overline\tau
:=e^{\|A\|_{\mathrm{op}}h_0}\|\Gamma\|_F,
\qquad
\sigma_\Delta\le\overline\sigma
:=e^{\|A\|_{\mathrm{op}}h_0}\|\Gamma\|_{\mathrm{op}}.
\]
Define
\[
C_{\mathrm{op}}:=\max\left\{1,
\frac{4M_\Phi M_x}{\lambda_\Phi}
\left(\overline\tau+\sqrt2\,\overline\sigma\right)\right\}.
\]
Since $\log(2/\delta)\le L_\delta$ and $L_\delta>1$,
\eqref{eq:md-concentration-radii} implies
\[
\varepsilon_{\mathrm{op},n,\Delta}\le
C_{\mathrm{op}}\sqrt{\frac{L_\delta}{n\Delta t}}.
\]
One valid choice of the constants in the theorem is
\begin{equation}
\label{eq:md-explicit-sample-constants}
\begin{split}
c_1&=\max\left\{
\frac{32M_\Phi^4}{\lambda_\Phi^2},
\frac{8M_\Phi^2}{3\lambda_\Phi}\right\},\qquad
c_2=\max\left\{1,\frac{C_{\mathrm{op}}^2}{{\varepsilon_0}^2}\right\},\\
C_{\mathrm{stat}}&=C_\star C_{\mathrm{op}}.
\end{split}
\end{equation}
Indeed, $n\ge c_1L_\delta$ makes each of the two terms in $r_G$
at most $\lambda_\Phi/4$, so $r_G\le\lambda_\Phi/2$.
The condition $n\Delta t\ge c_2L_\delta$ gives
$\varepsilon_{\mathrm{op},n,\Delta}\le\varepsilon_0$.

Proposition~\ref{prop:operator-concentration}
therefore supplies one event of probability at least $1-\delta$
on which $\widehat G$ is invertible and
\begin{equation}
\label{eq:md-operator-error-bound}
\|\widehat M-M_\Delta\|_{\mathrm{op}}
\le\varepsilon_{\mathrm{op},n,\Delta}
\le C_{\mathrm{op}}\sqrt{\frac{L_\delta}{n\Delta t}}
\le\varepsilon_0.
\end{equation}
Applying Theorem~\ref{thm:large-beta-operator} proves stability,
well-posedness, and quadratic convergence on this event.
Substituting \eqref{eq:md-operator-error-bound} into
its gain-error bound~\eqref{eq:md-deterministic-end-to-end} gives
\eqref{eq:md-large-beta-end-to-end} with the stated
$C_{\mathrm{stat}}$.
The covariance bounds and all chosen constants
are independent of $\beta,\Delta t,n,N,\delta$ over the stated range.
\end{proof}

\subsection{Continuous-time cost error}
\begin{corollary}[Continuous-time cost error]
\label{cor:md-cost}
Suppose the conditions of Theorem~\ref{cor:md-large-beta-end-to-end}
hold and the initial second moment $S_0=\E[X_0X_0^\top]$ is finite.
Define
$\Omega(\beta)=S_0+\Gamma\Gamma^\top/\beta$.
Let $J_{\mathrm{ct}}(K)$ be the expected discounted positive cost, and let
$\ell_{\mathrm{ct}}(\beta)$ bound the inverse Lyapunov operator of
$A-BK-\beta I/2$ on the gain ball $\mathcal U$ constructed in
Appendix~\ref{app:large-beta}, for which one may take
$\ell_{\mathrm{ct}}(\beta)=2/\beta$ by
\eqref{eq:large-beta-lyapunov-inverse} at $\Delta t=0$.
On the event of that theorem, simultaneously for every $N\ge1$,
\begin{equation}
\label{eq:md-cost-bound}
\begin{aligned}
0\le J_{\mathrm{ct}}(K_N)-J_{\mathrm{ct}}(K^\star_{\mathrm{ct}})\le\|\Omega(\beta)\|_F\ell_{\mathrm{ct}}(\beta)\|R\|_{\mathrm{op}}
\|K_N-K^\star_{\mathrm{ct}}\|_F^2.
\end{aligned}
\end{equation}
With $\Delta t\asymp n^{-1/3}$,
the cost gap is $O_p(n^{-2/3})$ after the iteration residual is negligible.
\end{corollary}

\begin{proof}[Proof of Corollary~\ref{cor:md-cost}]
Apply the completed-square identity in
\eqref{eq:md-quadratic-identity} to the true continuous-time matrices.
The cost is $J_{\mathrm{ct}}(K)=\operatorname{tr}(\Omega(\beta)\Pi_K^{\mathrm{ct}})$,
including the additive-noise value constant. The positive Lyapunov inverse
preserves positive semidefiniteness, proving nonnegativity; its norm bound
and the Frobenius Cauchy--Schwarz inequality give \eqref{eq:md-cost-bound}.
\end{proof}

\section{Bellman-Equation Baselines}
\label{app:baselines}
\subsection{BE-PI and its scalar population target}
\label{app:be-comparison}

The BE-PI baseline uses discount $\gamma=e^{-\beta \Delta t}$ and stage reward
$-\Delta t(qx^2+ru^2)$, where $\Delta t>0$. With the scalar features
$\Phi(x,u)=(1,x^2,xu,u^2)^\top$ and reward vector
$w=(0,-q,0,-r)^\top$, the BE residual is
$g_{K,i}^{\mathrm{BE}}=[\Phi_i-\gamma\Phi(x_i',-Kx_i')]/(\Delta t)$.
It forms $\widehat A_K^{\mathrm{BE}}=n^{-1}\sum_i\Phi_i(g_{K,i}^{\mathrm{BE}})^\top$
and solves $\widehat A_K^{\mathrm{BE}}\widehat\theta_K=\widehat G w$,
where $\widehat G=n^{-1}\sum_i\Phi_i\Phi_i^\top$.
Greedy improvement gives
$K^+=\widehat\theta_{K,xu}/(2\widehat\theta_{K,uu})$ when the action
coefficient is negative. The same batch is reused at every iteration.

\paragraph{Population target.}
Write the exact held-action transition as
\begin{equation}
\label{eq:be-transition}
x'=A_{\Delta t} x+B_{\Delta t} u+\epsilon,\quad
A_{\Delta t}=e^{A\Delta t},\quad B_{\Delta t}=\int_0^{\Delta t} e^{At}B\,dt,
\end{equation}
where $\epsilon$ is centered Gaussian with variance
$\sigma^2\int_0^{\Delta t} e^{2At}\,dt$.
For $\gamma(A_{\Delta t}-B_{\Delta t}K)^2<1$, the quadratic value coefficient is
$P_K=-\Delta t(q+rK^2)/[1-\gamma(A_{\Delta t}-B_{\Delta t}K)^2]$.
Substitution into the Bellman equation gives the population greedy map
$\mathcal G_{\mathrm{BE}}(K)=\gamma P_KA_{\Delta t}B_{\Delta t}/(-\Delta t\,r+\gamma P_KB_{\Delta t}^2)$.
Additive noise changes only the value constant and does not change this
map. Population BE-PI is the classical Newton--Hewer iteration
\citep{Hewer1971} for the discounted discrete-time algebraic Riccati equation with
left-endpoint stage costs $\Delta t\,q,\Delta t\,r$.
Its stabilizing optimal gain $K_{\mathrm{BE},\Delta t}^\star$ satisfies
\begin{equation}
\label{eq:be-root}
\gamma rA_{\Delta t}B_{\Delta t}(K_{\mathrm{BE},\Delta t}^\star)^2+
\bigl[r(1-\gamma A_{\Delta t}^2)+\gamma qB_{\Delta t}^2\bigr]K_{\mathrm{BE},\Delta t}^\star
-\gamma qA_{\Delta t}B_{\Delta t}=0.
\end{equation}
Thus its finite-step target differs from the surrogate CARE targeted by
PhiBE-PI. Related scalar discretization comparisons appear in
\citet[Theorem~4.4]{ZhuZhangZhang2025OptimalPhiBE}, with discount
$(1+\beta \Delta t)^{-1}$ rather than the exponential discount used here.

\subsection{Endpoint adaptation of Pang--Jiang policy iteration}
\label{app:pj-adaptation}

\citet{PangJiang2023Stationary} estimate an affine map from a quadratic
value matrix to the coefficients of the instantaneous cost plus its
generator, using continuous-time trajectory integrals.
We adapt this construction to the discounted
objective and endpoint transitions used here. The adaptation uses the
constant and homogeneous quadratic features $\Phi(x,u)$ from
Section~\ref{subsec:lqr-pi-algo}, omitting the original estimator's linear
terms, and replaces finite-time ODE evaluation with an exact solution of
the empirical evaluation equation.

Suppose the empirical Gram matrix $\widehat G$ is
invertible. For a state function $f$, define the successor regression
\begin{equation}
\label{eq:pj-successor-regression}
(\widehat{\mathcal P}f)(x,u)
=\Phi(x,u)^\top\widehat G^{-1}\frac1n\sum_i\Phi_i f(X_i').
\end{equation}
Regression reproduces current-state quadratics exactly because they
belong to the feature space. For a candidate value $V(x)=c+x^\top Ux$,
the endpoint adaptation fits $R_i+[V(X_i')-V(X_i)]/(\Delta t)$, giving
\begin{equation}
\label{eq:pj-endpoint-generator}
\widehat{\mathcal A}_V(x,u)=\mathcal R(x,u)
+\frac{(\widehat{\mathcal P}V)(x,u)-V(x)}{\Delta t}.
\end{equation}
Policy evaluation requires $\widehat{\mathcal A}_V(x,-Kx)=\beta V(x)$,
or equivalently
\begin{equation}
\label{eq:pj-value-evaluation}
(1+\beta \Delta t)V(x)=\Delta t\mathcal R(x,-Kx)+(\widehat{\mathcal P}V)(x,-Kx).
\end{equation}
This is linear in the value coefficients. Policy improvement maximizes
$\widehat{\mathcal A}_V(x,u)$ over $u$, provided its quadratic action
block is negative definite.

This adaptation also corresponds to the modified
Optimal-BE formulation of
\citet[Eq.~(38)]{ZhuZhangZhang2025OptimalPhiBE}, with discount
$\gamma=(1+\beta\Delta t)^{-1}$ and stage reward
$\gamma\Delta t\,\mathcal R$.

\begin{proposition}[Fixed-batch equivalence]
\label{prop:pj-be-equivalence}
Suppose the empirical Gram matrix is invertible. Use the same batch,
feature space, and initial gain for PJ-PI (adapted) and BE-PI with
stage reward $\Delta t\mathcal R$ and discount $\gamma=(1+\beta \Delta t)^{-1}$.
If policy evaluation is solved exactly, the two methods produce
identical gain updates whenever the evaluation solutions are unique
and the action blocks are negative definite. Consequently, their gain
sequences coincide for as long as these conditions hold.
\end{proposition}

\begin{proof}
Put $\alpha=1+\beta \Delta t$, $\gamma=\alpha^{-1}$, and define
$Q_B(x,u)=\Delta t\widehat{\mathcal A}_V(x,u)+V(x)$.
The adapted evaluation equation gives
$V_B(x):=Q_B(x,-Kx)=\alpha V(x)$. Thus
\[
Q_B=\Delta t\mathcal R+\widehat{\mathcal P}V
=\Delta t\mathcal R+\gamma\widehat{\mathcal P}V_B,
\]
which is the empirical BE evaluation equation; multiplication by
$\widehat G$ recovers its LSTD moment system.
Conversely, a BE solution gives an adapted evaluation solution by
$V=V_B/\alpha$. Since $\Delta t>0$ and $V(x)$ is independent of $u$,
$Q_B$ and $\widehat{\mathcal A}_V$ have the same maximizing action.
Induction from the common initial gain proves the claim.
\end{proof}

\section{Experimental Details and Supplementary Results}
\label{app:additional-experiments}

\subsection{Motivating example: temporal resolution for learning}
\label{app:motivation-setup}
Figure~\ref{fig:learning-temporal-resolution} uses the scalar system
$dX_t=(-2X_t+u_t)\,dt+dW_t$, with $Q=512$, $R=1$, and $\beta=5$.
Each trial starts from $X_0\sim\mathcal N(0,1/4)$ and contains $64$
consecutive constant-action segments. At the start of each segment, an
action is drawn independently and uniformly from $[-8,8]$ and held for
$8192$ fine steps. States are carried continuously between segments.
We simulate the exact Gaussian held-action transitions at observation
interval $\delta=0.002$, giving $N=524{,}288$ fine-step transitions per trial.

For each
\[
k\in\{1,4,8,12,16,20,24,32,64,128,256,384,512,768,1024\},
\]
we construct all
pairs $(X_{j\delta},u_j,X_{(j+k)\delta})$ whose endpoints belong to the
same constant-action segment. These pairs overlap, and their number is
$n_k=64(8193-k)$, retaining between $87.5\%$ and $100\%$ of $N$.
The same data are used by both methods and at every temporal resolution
$\Delta t=k\delta$. Trials are independent, while pairs within a trial
share states and process noise.

Both methods use the features $\Phi(x,u)=(1,x^2,xu,u^2)^\top$,
initial gain $K_0=1$, and unregularized policy-evaluation solves.
Each method reuses its fixed batch throughout policy iteration and stops
when $|K_{\ell+1}-K_\ell|\le10^{-11}(1+|K_\ell|)$, with at most
$60$ updates. The discrete-time Bellman-equation policy-iteration
baseline (BE-PI) uses discount $e^{-\beta\Delta t}$ and stage reward
$-\Delta t(512x^2+u^2)$, as in Appendix~\ref{app:be-comparison}.
Deployment always uses $u=-KX$ with the state observed and the action
updated every $\delta$, independently of $\Delta t$.

A fit fails if a policy update is nonfinite, the critic loses strict
concavity in the action, or the iteration does not converge within
$60$ updates. Both methods retain the supplied stable gain $K_0$ on a
failed fit. This rule uses only observable training failures; no trial
is discarded and true-model performance is not used to choose a fallback.

For evaluation, let
\[
J_\delta(K)=\mathbb E\!\left[\int_0^\infty e^{-5t}
\bigl(512X_t^2+u_t^2\bigr)\,dt\right],
\qquad u_t=-KX_{\delta\lfloor t/\delta\rfloor},
\]
with evaluation initial state $X_0\sim\mathcal N(0,1)$.
The reference $J_\delta^\star=\min_K J_\delta(K)\approx22.35459$
is the optimal cost among linear feedback gains deployed at the same
interval $\delta$. It is computed from the discounted discrete-time
Riccati equation with the exact within-interval integrated cost, including
the process-noise contribution. The figure reports mean relative cost errors
\[
e(\widehat K)=100\,\frac{J_\delta(\widehat K)-J_\delta^\star}
{J_\delta^\star}\quad(\%)
\]
over $1000$ paired validation trials.
Shaded bands are pointwise $95\%$ confidence intervals for the mean,
computed as $\bar e\pm1.96s_e/\sqrt{1000}$, where $s_e$ is the
sample standard deviation across trials. The main axes are logarithmic;
the inset compares the highlighted intervals on a linear vertical scale
starting at zero.

Independent pilot data select $\Delta t=0.032$ for BE-PI and
$\Delta t=1.536$ for our algorithm before generating the validation trials.
At these fixed intervals, the mean relative cost errors are $7.83\%$
and $2.67\%$, respectively, a reduction of $65.9\%$.
There are $8$ failed fits for BE-PI and none for our algorithm, and all
deployed gains at these intervals are stable. At the finest training
interval, $\Delta t=\delta$, the failure counts are $167$ and $190$,
respectively. Thus both curves include a training-failure penalty at fine
resolutions and exhibit a training-resolution tradeoff, despite the fixed
deployment interval. The highlighted intervals were selected before
validation and are not presented as validation-grid minima.

\subsection{Data generation and implementation}
\label{app:experiment-implementation}
The experiments in Section~\ref{subsec:lqr-numerics-main} use
independent state--action starting pairs and exact
Gaussian held-action transitions. Within each trial, the three methods
share the same batch and initial gain, and reuse the batch throughout
policy iteration. True dynamics are used only for data generation and
performance evaluation. Policy evaluation uses a constant and all
homogeneous quadratic features, with exact unregularized solves.

Policy iteration stops when
$\|K_{k+1}-K_k\|_F\le\varepsilon(1+\|K_k\|_F)$, with at most
30 updates. The convergence experiment below instead
retains a fixed number of updates to display the numerical error floor.
We use $\varepsilon=10^{-10}$ for the scalar experiments in
Figures~\ref{fig:beta-dt-control} and~\ref{fig:deployment-ratio-control},
and $10^{-9}$ for the remaining experiments. Runs terminate upon a
failed solve, an evaluation condition number above $10^{12}$, or loss
of negative definiteness of the action block. All runs reported below
converge and have finite discounted cost.

The scalar experiments use $A=0.4$,
$B=\Gamma=Q=R=1$, $K_0=1.5$, and starting state--action pairs uniform
on $[-2,2]^2$. The scalar features are
$\Phi(x,u)=(1,x^2,xu,u^2)^\top$. Reference gains are obtained from the
surrogate CARE for PhiBE-PI and the corresponding discounted DAREs
for the two baselines.

Figure~\ref{fig:beta-dt-control} uses 30 paired trials;
Figure~\ref{fig:deployment-ratio-control}(a) uses 200 paired trials.
Sample-size sweeps use nested batches; discount and observation-interval
sweeps share starting pairs and standard Gaussian innovations within
each trial. Means and pointwise 95\% confidence
intervals are computed using 1000 bootstrap resamples of complete
trial paths.

\subsection{Convergence to the empirical fixed point}
\label{app:convergence-experiment}
\begin{figure}[!htbp]
\centering
\includegraphics[width=0.9\linewidth]{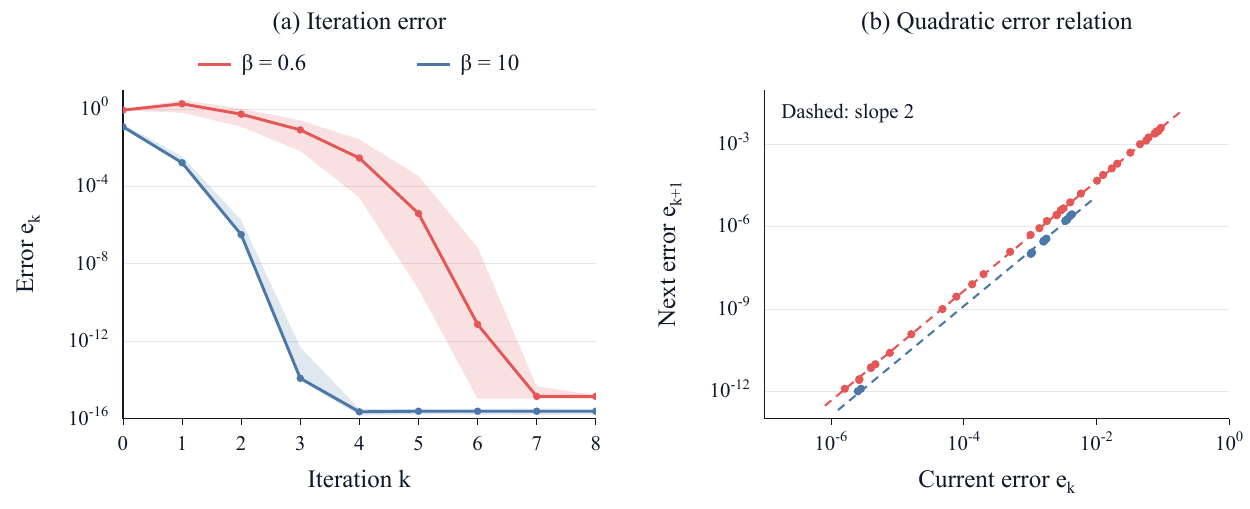}
\caption{Convergence to the empirical fixed point.
(a) Iteration errors relative to
high-precision empirical fixed points, with medians and interquartile
ranges over 15 batches per discount. (b) Consecutive errors in the
local regime, with slope-two reference lines. All 30 runs converge.}
\label{fig:pi-convergence-diagnostics}
\end{figure}

Figure~\ref{fig:pi-convergence-diagnostics} uses the same five normal systems as
Section~\ref{app:multid-deployment}, with $Q=I_4$, $R=I_2$,
$\Gamma=0.35I_4$, and $K_0=0$. For each system we generate three
independent batches of $n=32768$ exact held-action transitions at
$\Delta t=0.1$, starting uniformly on $[-1,1]^6$. The two discounts
$\beta\in\{0.6,10\}$ share these batches. Every run reuses its batch
for 12 unregularized policy-iteration updates in double precision;
Figure~\ref{fig:pi-convergence-diagnostics}(a) displays iterations $0$ through $8$.
All 30 runs converge, with positive definite cost-critic action blocks
and discounted-stabilizing gains for both the true and surrogate
systems. The cost-critic convention reverses the sign of the reward
critic used in the main text.

We measure iteration error as
$e_k=\|K_k-\widehat K_\Delta\|_F$, using an empirical fixed point
refined from each stored evaluation system with 70-digit arithmetic.
We check its residual with 90-digit arithmetic; the largest is below
$3\times10^{-69}$. The plotted iterates use double precision, and the
plateaus near $10^{-15}$--$10^{-16}$ reflect numerical roundoff.

Figure~\ref{fig:pi-convergence-diagnostics} gives the full diagnostics.
Panel~(a) reports medians and interquartile ranges over all 15 batches
per discount. Panel~(b) includes every
consecutive pair satisfying $e_k<0.1$ and
$10^{-12}<e_{k+1}<e_k$, restricting attention to the local regime
above the numerical error floor. The dashed lines have slope two,
with their intercepts set by the geometric mean of
$e_{k+1}/e_k^2$ for each discount. For eligible triples with
$e_{k-1}>e_k$, the observed orders
$\log(e_{k+1}/e_k)/\log(e_k/e_{k-1})$ have medians $1.989$ and
$1.988$ for $\beta=0.6$ and $10$, respectively.
The update tolerance $10^{-12}(1+\|K_k\|_F)$ is first met after
6--9 updates for $\beta=0.6$ and 4--5 updates for $\beta=10$.

\subsection{Gain accuracy and sample-size scaling}
\label{app:scaling-details}
Figure~\ref{fig:beta-dt-control}(a) fixes $n=10^8$ and uses
$\Delta t\in\{0.005,0.02,0.1\}$, with 31 logarithmically spaced
discounts from $0.1$ to $100$, supplemented by $\{0.5,1,5,10,20\}$.
Panel~(b) uses $n=10^8$, $\beta\in\{1,5,10\}$, and
21 logarithmically spaced intervals from $0.001$ to $0.5$,
supplemented by $\{0.005,0.01,0.02,0.1\}$. Panel~(c) uses
$\beta\in\{1,10\}$ and nine logarithmically spaced sample sizes
from $10^4$ to $10^8$, rounded to integers, at fixed $\Delta t=0.02$.

Table~\ref{tab:scaling-slopes} reports log--log
slopes fitted to mean errors over the nine sample sizes and 30 paired
trials in Figure~\ref{fig:beta-dt-control}(c). The fixed-interval sweep
$\Delta t=0.02$ measures finite-sample gain error relative to each
method's population target. The balanced-interval sweep
$\Delta t=n^{-1/3}$ uses the same sampled prefixes and measures total
gain error relative to $K^\star_{\mathrm{ct}}$ and relative continuous-time
cost gap. Both gain errors use absolute differences.

\begin{table}[!htbp]
\centering
\caption{Empirical log--log slopes over nine sample sizes
from $10^4$ to $10^8$. Finite-sample gain error uses $\Delta t=0.02$;
total gain error and relative continuous-time cost gap use the separate
balanced-interval sweep $\Delta t=n^{-1/3}$.
Reference gives the theoretical exponent for PhiBE-PI.}
\label{tab:scaling-slopes}
\begingroup\small
\setlength{\tabcolsep}{4pt}
\renewcommand{\arraystretch}{1.25}
\begin{tabular}{llcccc}
\toprule
Metric & $\beta$ & Reference & PhiBE-PI & BE-PI & PJ-PI (adapted)\\
\midrule
Finite-sample gain error & 1 & $-1/2$ & $-0.504$ & $-0.504$ & $-0.504$\\
 & 10 & $-1/2$ & $-0.521$ & $-0.520$ & $-0.520$\\
Total gain error & 1 & $-1/3$ & $-0.337$ & $-0.351$ & $-0.346$\\
 & 10 & $-1/3$ & $-0.349$ & $-0.330$ & $-0.341$\\
Relative CT cost gap & 1 & $-2/3$ & $-0.685$ & $-0.718$ & $-0.709$\\
 & 10 & $-2/3$ & $-0.695$ & $-0.659$ & $-0.688$\\
\bottomrule
\end{tabular}
\endgroup
\end{table}

The slopes are consistent with the $n^{-1/2}$ finite-sample scale and
the $n^{-1/3}$ gain and $n^{-2/3}$ cost-gap bounds for PhiBE-PI.
Both baselines exhibit similar sample-size scaling in this sweep.

\FloatBarrier
\subsection{Training and deployment frequencies}
\label{app:multid-deployment}
Figure~\ref{fig:deployment-ratio-control} fixes the gains learned at
$\Delta_{\mathrm{tr}}=0.5$ and varies the deployment interval
$\Delta_{\mathrm{ts}}=\Delta_{\mathrm{tr}}/\rho$. Performance is measured
against the optimal held-action policy at each deployment interval.
Costs integrate the discounted running cost over each action-hold
interval, including the state--action cross term and process-noise
contribution. Both baselines use left-endpoint stage costs during
training, and PJ-PI (adapted) also modifies the one-step discount,
so their gaps need not vanish at $\rho=1$.

Panel~(a) uses the scalar system with $\beta=0.5$, $n=2^{18}$,
initial second moment $1$, and $\rho\in\{1,2,5,10,20,50,100\}$.
Interval costs are computed using 64-point Gauss--Legendre quadrature.

Panels~(b,c) use four-state, two-input systems with $Q=I_4$, $R=I_2$,
and $K_0=0$. The two families are
constructed as follows:
\begin{itemize}
\item \textbf{Normal systems:} $A=UDU^\top$ and $B=UB_0$, where
$U$ is the orthogonal factor of a Gaussian matrix, $D$ has four
equally spaced diagonal entries from $-0.8$ to $0.15$, and $B_0$
has independent $N(0,1/4)$ entries. We use $\beta=0.6$ and
$\Gamma=0.35I_4$ during both training and deployment.
\item \textbf{Mass--spring systems:}
$A=\left[\begin{smallmatrix}0&I\\-0.1L&-8I\end{smallmatrix}\right]$
and $B=\left[\begin{smallmatrix}0\\I\end{smallmatrix}\right]$,
where $L_{ii}=k_i+k_{i+1}+0.2$,
$L_{i,i+1}=L_{i+1,i}=-k_{i+1}$, and the $k_i$ are independent
and uniform on $[0.6,1.4]$. We use $\beta=0.1$ and
$\Gamma=0.35I_4$; the resulting systems are overdamped.
\end{itemize}

Each family contains five systems. The normal systems use three paired
batches of $10^8$ transitions per system; the mass--spring systems use
six paired batches of $10^8$ transitions per system, giving 30 trials.
Starting state--action pairs are uniform on $[-1,1]^6$.
Transitions are generated from the exact held-action linear dynamics,
including the Gaussian transition covariance for both families.
All three methods reuse the same batch throughout policy iteration.
Deployment intervals are
$\{0.01,0.02,0.05,0.1,0.2,0.5\}$, with initial covariance $I_4$.
Panels show medians and interquartile ranges across systems after
averaging relative cost gaps over batches within each system.
Section~\ref{app:lightly-damped-mass-spring}
examines the effect of first-order bias in a lightly damped system.

\subsection{A lightly damped mass--spring example}
\label{app:lightly-damped-mass-spring}
We retain the unscaled stiffness matrix $L$ defined above and use
$A=\left[\begin{smallmatrix}0&I\\-L&-0.2I\end{smallmatrix}\right]$
and $B=\left[\begin{smallmatrix}0\\I\end{smallmatrix}\right]$,
with $Q=I_4$, $R=I_2$, $\beta=0.6$, and $K_0=0$.
To isolate finite-step bias, this example uses noiseless transitions
($\Gamma=0$), five systems, and three paired batches of $n=8192$
transitions per system. Starting pairs are uniform on $[-1,1]^6$.
The training interval, deployment intervals, initial covariance, and
cost reference are as in Section~\ref{app:multid-deployment}.

\begin{figure}[!t]
\centering
\includegraphics[width=0.58\linewidth]{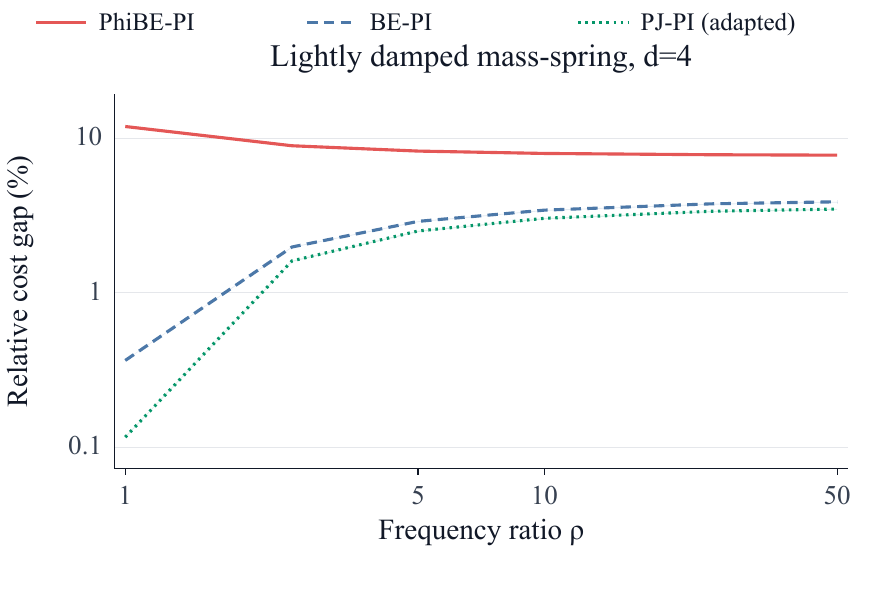}
\caption{A lightly damped mass--spring example in which
first-order PhiBE-PI has a larger deployment cost gap than both baselines.
Training uses $\Delta_{\mathrm{tr}}=0.5$, $\beta=0.6$, $\Gamma=0$,
and $n=8192$. Curves and bands show medians and interquartile ranges
across five systems after averaging three paired batches per system.
The reference is the optimal held-action policy at each deployment interval.}
\label{fig:lightly-damped-mass-spring}
\end{figure}

Figure~\ref{fig:lightly-damped-mass-spring}
shows a larger deployment cost gap for PhiBE-PI: at $\rho=50$, the
median gaps are $7.79\%$, $3.88\%$, and $3.49\%$ for PhiBE-PI,
BE-PI, and PJ-PI (adapted), respectively. All 45 fits converge and agree
with their population CARE or DARE targets to numerical precision.

The first-order population approximation perturbs both drift matrices:
\[
\bar A_\Delta=A+\frac{\Delta t}{2}A^2+O(\Delta t^2),
\qquad
\bar B_\Delta=B+\frac{\Delta t}{2}AB+O(\Delta t^2).
\]
In lightly damped oscillatory modes, the correction to $A$ can introduce
substantial artificial decay. For one of the tested systems, eigenvalues
$-0.100\pm1.757\mathrm{i}$ become $-0.786\pm1.464\mathrm{i}$ in
$\bar A_\Delta$ at $\Delta t=0.5$.
The bias in $\bar B_\Delta$ also affects the learned control.
Faster deployment of the fixed gain cannot remove these training-time
biases. This example and the overdamped setting in
Figure~\ref{fig:deployment-ratio-control}(c) motivate systematic
benchmarking of the regimes in which PhiBE is effective.

\subsection{Noisy multidimensional experiments}
\label{app:paired-noisy-comparison}
We also compare the methods on the normal systems above and lightly
damped mass--spring systems with process noise, using the scalar system
as a reference. In this supplementary comparison, the mass--spring drift
is $A=\left[\begin{smallmatrix}0&I\\-L&-0.2I\end{smallmatrix}\right]$,
with the unscaled $L$ defined above. All experiments use
$n=32768$, $\Delta t=0.1$, diffusion $\Gamma=0.35I$, $Q=I$,
$R=I$, and initial covariance $I$. The scalar system uses $\beta=0.5$,
$K_0=1.5$, and starting pairs uniform on $[-2,2]^2$; the four-state
systems use $\beta=0.6$, $K_0=0$, and starting pairs uniform on
$[-1,1]^6$.

We use 100 paired trials for the scalar system and five systems with
20 paired batches each for each four-state family.
Table~\ref{tab:paired-noisy-comparison} reports medians over all
100 runs per family.

\begin{table}[!htbp]
\centering
\caption{Noisy scalar and four-state systems.
Gain error is $100\|K-K^\star_{\mathrm{ct}}\|_F/\|K^\star_{\mathrm{ct}}\|_F$.
Relative CT cost gap includes process noise. All entries are percentages.
Bold entries indicate the lowest error within each system family.}
\label{tab:paired-noisy-comparison}
\begingroup\small
\begin{tabular}{llrr}
\toprule
Family & Method & Gain error & CT cost gap\\
\midrule
Scalar & PhiBE-PI & \textbf{0.781} & \textbf{0.0035}\\
 & BE-PI & 4.533 & 0.1245\\
 & PJ-PI (adapted) & 3.958 & 0.0942\\
\midrule
Normal & PhiBE-PI & \textbf{7.363} & \textbf{0.1054}\\
 & BE-PI & 8.681 & 0.1232\\
 & PJ-PI (adapted) & 8.488 & 0.1267\\
\midrule
Mass--spring & PhiBE-PI & 16.855 & 0.8467\\
 & BE-PI & 15.603 & 0.3088\\
 & PJ-PI (adapted) & \textbf{14.971} & \textbf{0.2801}\\
\bottomrule
\end{tabular}
\endgroup
\end{table}

PhiBE-PI has the smallest median gain error and cost gap for the scalar
and normal systems, whereas PJ-PI (adapted) has the smallest medians
for the mass--spring systems.

\end{document}